\documentclass[12pt]{article} 

\usepackage{epsfig,amsmath,amssymb,graphics,verbatim,amsfonts,subfigure,psfrag}
\usepackage[letterpaper,margin=0.9in]{geometry}
\usepackage{algorithmic,color}
\usepackage[lined,boxed,ruled,commentsnumbered]{algorithm2e}
\usepackage{titlesec}

\def\txtb{{\textnormal{b}}}

\def\txtd{{\textnormal{d}}}
\def\txte{{\textnormal{e}}}
\def\txtf{{\textnormal{f}}}

\def\R{\mathbb{R}}
\def\C{\mathbb{C}}
\def\N{\mathbb{N}}
\def\P{\mathbb{P}}
\def\E{\mathbb{E}}
\def\I{\infty}

\newcommand{\be}{\begin{equation}}
\newcommand{\ee}{\end{equation}}
\newcommand{\bea}{\begin{eqnarray}}
\newcommand{\eea}{\end{eqnarray}}
\newcommand{\beann}{\begin{eqnarray*}}
\newcommand{\eeann}{\end{eqnarray*}}
\newcommand{\benn}{\begin{equation*}}
\newcommand{\eenn}{\end{equation*}}

\def\ra{\rightarrow}
\def\I{\infty}

\newcommand{\cA}{{\mathcal A}}  
\newcommand{\cB}{{\mathcal B}}  
\newcommand{\cC}{{\mathcal C}}  
\newcommand{\cD}{{\mathcal D}}  
\newcommand{\cF}{{\mathcal F}}  
\newcommand{\cG}{{\mathcal G}}  
\newcommand{\cH}{{\mathcal H}}  
\newcommand{\cI}{{\mathcal I}}  
\newcommand{\cK}{{\mathcal K}}  
\newcommand{\cL}{{\mathcal L}}  
\newcommand{\cM}{{\mathcal M}}  
\newcommand{\cO}{{\mathcal O}}  
\newcommand{\cU}{{\mathcal U}}  
\newcommand{\cV}{{\mathcal V}}  

\begin{document}

\titleformat*{\section}{\large\bfseries}

\author{Christian Kuehn\thanks{Institute for Analysis and 
Scientific Computing, Vienna University of Technology, 1040 Vienna, Austria}}

\title{Numerical Continuation and SPDE Stability\\ for the 2D Cubic-Quintic Allen-Cahn Equation}

\maketitle

\begin{abstract}
We study the Allen-Cahn equation with a cubic-quintic nonlinear term and a $Q$-trace-class 
stochastic forcing in two spatial dimensions. This stochastic partial
differential equation (SPDE) is used as a test case to understand, how numerical 
continuation methods can be carried over to the SPDE setting. First, we compute the deterministic
bifurcation diagram for the PDE, {i.e.}~without stochastic forcing. In this case, two locally asymptotically
stable steady state solution branches exist upon variation of the linear damping term. Then
we consider the Lyapunov operator equation for the locally linearized system around steady
states for the SPDE. We discretize the full SPDE using a combination of finite-differences
and spectral noise approximation obtaining a finite-dimensional system of stochastic
ordinary differential equations (SODEs). The large system of SODEs is used to approximate
the Lyapunov operator equation via covariance matrices. The covariance matrices are numerically
continued along the two bifurcation branches. We show that we can quantify the stochastic
fluctuations along the branches. We also demonstrate scaling laws near branch and fold bifurcation points.
Furthermore, we perform computational tests to show that, even with a sub-optimal computational
setup, we can quantify the subexponential-timescale fluctuations near the deterministic
steady states upon stochastic forcing on a standard desktop computer setup. Hence, the
proposed method for numerical continuation of SPDEs has the potential to allow for rapid parametric 
uncertainty quantification of spatio-temporal stochastic systems.
\end{abstract}

{\bf Keywords:} Allen-Cahn, SPDE, cubic-quintic Ginzburg-Landau, Lyapunov equation, numerical 
continuation, spectral noise approximation, uncertainty quantification, finite differences, 
iterative linear solvers, critical transitions, bifurcation diagram.\\

\section{Introduction}
\label{sec:intro}

Spatio-temporal systems with stochastic forcing \emph{and} parameter dependence arise extremely
frequently in practical applications. Stochastic partial differential equations (SPDEs) have
emerged as one of the main examples for spatio-temporal stochastic systems, 
{e.g.}~in ecological invasion waves \cite{ElworthyZhaoGaines,KuehnFKPP,MuellerSowers}, 
voter models \cite{MuellerTribe1,CoxDurrettPerkins1}, 
surface and thin film growth models \cite{BloemkerRomito,BloemkerHairer1}, 
statistical mechanics \cite{KardarParisiZhang,Hairer2},
fluid dynamics \cite{HohenbergSwift,Bensoussan1},
neuroscience \cite{Tuckwell2,SauerStannat}, 
and many general pattern formation problems  
\cite{Bloemker,GarciaOjalvoSancho}. In this paper, we study as a prototypical model a version of the 
Allen-Cahn SPDE, which can formally be written as
\be
\label{eq:AC_intro}
\frac{\partial u}{\partial t}=\Delta u-F(u;\mu)+G(u)\eta,
\ee
where $u=u(z,t)\in\R$, $z\in\R^2$, $\mu\in\R^p$ are parameters, $F,G$ are suitable maps, $\Delta$ is 
the Laplacian and $\eta$ is a stochastic process. Our main interest is to use \eqref{eq:AC_intro} as 
a test problem to illustrate a different approach to quantify stochastic fluctuations using numerical 
continuation methods. A more precise formulation of \eqref{eq:AC_intro} and the detailed reasons
for choosing a certain Allen-Cahn-type SPDE are discussed below.\medskip

Currently, a standard approach to study the dynamics of \eqref{eq:AC_intro} numerically is centered 
around Monte-Carlo (MC) algorithms and their variations. The strategy is to use forward integration 
of the SPDE. This requires three different components:

\begin{enumerate}
 \item[(MC1)] Sample over a sufficiently large number of initial conditions $u(z,0)$ and determine 
 which initial conditions lead to similar dynamics.
 \item[(MC2)] Sample over a sufficiently large number of noise realizations for the process $\eta$ 
 and compute the required averaged quantities of interest over several sample paths.
 \item[(MC3)] Sample over the parameters $\mu \in\R^p$ and plot the bifurcation diagram at each point 
 by repeating steps (MC1)-(MC2). 
\end{enumerate} 

There are several strategies that partially alleviate some of the computational time required for 
carrying out (MC1)-(MC3), such as quasi-Monte-Carlo methods \cite{Caflisch,Niederreiter} based on 
discrepancy theory or multilevel-Monte-Carlo \cite{Heinrich,Giles} based upon a coarse-/fine-sampling 
hierarchy. However, there are some drawbacks when following an MC-based strategy. Due to the three
sampling spaces (initial conditions, noise, parameters), the dimension of the computational problem
can become quite restrictive, particularly for the case of SPDEs. Furthermore, inspection of 
individual sample paths is often necessary to study dynamical objects such as bifurcations, 
metstable paths, or locally invariant manifolds. This can make MC methods in practical applications 
surprisingly labor-intensive. Even for deterministic differential equations, such as ordinary differential 
equations (ODEs) or partial differential equations (PDEs), this issue has been known for some time as 
discussed in \cite[p.1-2]{TuckermanBarkley}.\medskip

Numerical continuation for ODEs and PDEs provides a systematic alternative to study dynamical 
phenomena. The idea is to track dynamical objects based upon parameter variation by a predictor-corrector 
scheme \cite{Keller,AllgowerGeorg,Govaerts}. This 
approach is used successfully in a wide variety of applications to study the parameter depedence of models 
\cite{KrauskopfOsingaGalan-Vioque,KuepperSeydelTroger,Seydel}. In combination with singularity detection 
algorithms, a common application of numerical continuation methods is to compute bifurcation diagrams 
\cite{Doedel1,Kuznetsov}. For ODE bifurcation analysis, excellent software packages are available 
such as \texttt{AUTO} \cite{Doedel_AUTO2007} or \texttt{MatCont} \cite{DhoogeGovaertsKuznetsov}. Although
there are many other applications of numerical continuation methods (see {e.g.}~the short summary
in the introduction to \cite{KuehnEllipticCont}), there is relatively little work, which tries to
apply these ideas to stochastic problems \cite{BarkleyKevrekidisStuart,KuehnSDEcont1}.\medskip

A possible strategy to make numerical continuation methods useful for stochastic ordinary differential
equations (SODEs) suggested in \cite{KuehnSDEcont1} is composed of the following main steps:

\begin{enumerate}
 \item[(N1)] Continue the deterministic problem without forcing with techniques for ODEs. For example, 
 one might want to compute the steady states and their bifurcations depending upon parameters.
 \item[(N2)] Use stochastic analysis on the level of the SODE to derive deterministic algebraic and/or 
 differential equations for important stochastic quantities. For example, one may derive an ODE
 for the evolution of local stochastic fluctuations around a steady state for the linearized system. Another
 example are algebraic equations for the pre-factor and exponent in Kramers' law governing metastable
 gradient systems.
 \item[(N3)] Use the deterministic equations from (N2) as \emph{additional} parts of the numerical 
 continuation problem. For example, one may use the numerical methods in (N1) also for the ODE derived
 in (N2).
\end{enumerate}

There are several advantages of the approach (N1)-(N3). The main advantage is that one works directly
with the dynamical objects of interest. For example, if one is interested in the local stochastic 
fluctuations on a certain time scale depending upon parameters, and one finds 
an algebraic equation to describe this behaviour in step (N2), then one just focuses the numerical 
computations upon this task. This generally yields substantial computational- and labour-time
savings. Furthermore, the step (N1) can be performed independently, if necessary. This allows for 
post-processing of models, which have been studied on the deterministic level already via numerical
continuation. The disadvantage of (N1)-(N3) is that one actually has to invest a lot of analytical
work in (N2) to determine the right equations and in (N3) to determine good discretizations for
these equations. However, this effort is mitigated by the fact that this work in (N2)-(N3) has only 
to be performed \emph{once}. The situation is similar to the ODE case. One could try to compute
bifurcation diagrams for ODEs via MC methods (also called ``brute-force'' bifurcation diagrams 
\cite{LinaroChampneysDesrochesStorace}). However, one may also invest quite some effort to find the correct
algebraic equations for numerical detection and continuation of bifurcations \cite{Kuznetsov} and then
automate the strategy, which leads to very fast and flexible software tools 
\cite{DhoogeGovaertsKuznetsov,Doedel_AUTO2007}.\medskip

In this paper, we give a first \emph{proof-of-concept} that the ideas outlined in (N1)-(N3) have strong
potential to work for SPDEs. In particular, we study the Allen-Cahn SPDE \eqref{eq:AC_intro} with a 
polynomial nonlinearity 
\be
\label{eq:nonlin_main}
F(u;\mu)=4(\mu u+u^3-u^5), 
\ee
where $\mu\in\R$ is the main bifurcation parameter. Next, we briefly summarize the main results and 
outline the structure of the paper. In Section \ref{sec:AC}, some background for 
the deterministic Allen-Cahn PDE is provided. We also motivate the choice \eqref{eq:nonlin_main} 
of the nonlinearity in more detail. Section \ref{sec:numcont_basics} very briefly reviews numerical 
continuation methods for PDEs to fix some notation. In Section \ref{sec:det_cont}, we use the toolbox
\texttt{pde2path} \cite{UeckerWetzelRademacher} to compute the main solution branches upon variation
of $\mu$ for the Allen-Cahn PDE on an asymmetric rectangular domain. This computation is well-known
and yields branch points from a homogeneous solution $u\equiv 0$ as well as fold points on the non-trivial 
branches. In particular, we are interested in those parts of the homogeneous branch $\Gamma_0$ and 
of the first non-trivial branch $\Gamma_1$, which are locally asymptotically stable for the Allen-Cahn PDE.

In Section \ref{sec:AC_SPDE_intro}, we discuss the Allen-Cahn SPDE \eqref{eq:AC_intro}. We use a 
$Q$-trace-class Wiener process to define the stochastic forcing term $\eta$ following the standard
Hilbert space evolution equation approach to SPDEs. The main goal for the SPDE \eqref{eq:AC_intro}
is to find an efficient numerical method to quantify the stochastic fluctuations depending upon 
$\mu$ and $G$ in phase space regions near deterministic locally asymptotically stable solutions. 
In this paper,
we just focus on the regime of sub-exponential time scales, {i.e.}, before any large deviation
effects have occured. To capture the fluctuations, we consider in Section \ref{sec:abs_Lyapunov}
a suitable Lyapunov operator equation, which captures the covariance of the process. Instead of 
just discretizing the operator equation, we 
take a more general approach and start from the full Allen-Cahn SPDE in Section \ref{sec:fd}.
The Laplacian is discretized using a finite-difference scheme, the noise term is approximated
by using truncation of its spectral representation and the resulting equation is projected onto
piecewise constant basis functions. This approach yields a large set of SODEs. We already note 
that any other spatial discretization technique would also work in the approach (N1)-(N3). 
In Section 
\ref{sec:LyaSODEs}, we discuss the covariance matrix of a suitable linearized SODE and the 
associated Lyapunov matrix equation. The covariance matrix for the SODEs is the required 
approximation for the solution of the full Lyapunov operator equation.

The finite-dimensional Lyapunov equation is solved numerically using iterative solvers
in a predictor-corrector framework along parts of the bifurcation branches $\Gamma_0$ and 
$\Gamma_1$ in Section \ref{sec:contSPDE_Lya}. The main numerical results are then presented
and discussed in Section \ref{ssec:lower_res} for the homogeneous branch and in Section 
\ref{ssec:upper_res} for the non-trivial branch. We summarize the main conclusions in a 
non-technical way:

\begin{enumerate}
 \item[(C1)] The computation of the large covariance matrices can be carried out along a typical
bifurcation branch for the Allen-Cahn SPDE within a few minutes on a current standard desktop computer
setup. This holds despite the fact that we have not used specialized numerical tools and 
only used a very simple predictor-corrector scheme. Hence, there is strong
potential of the general approach (N1)-(N3) to allow for rapid uncertainty quantification analysis.
Furthermore, a comparison between different large-scale linear solvers has been carried out for
this problem.
 \item[(C2)] The stochastic fluctuations increase as branch and fold bifurcation points are approached along 
 a stable branch as 
expected from the theory of critical transitions. The inverse linear scaling for the branch point and
the square-root scaling for the fold point were observed in accordance with theoretical predictions.
 \item[(C3)] A numerical comparison between different norms to measure fluctuations, noise truncation
levels and various additive as well as multiplicative noise terms has been performed. This shows
the flexibility of the approach to characterize the influence of different noise terms on the local 
dynamics near parameterized steady state branches.  
\end{enumerate}

Before proceeding to the main part of the manuscript, we remark that there is a very large
number of open problems that could be treated using the philosophy (N1)-(N3) in the context
of SPDEs. Here we have just restricted ourselves to a \emph{proof-of-concept} to present 
the main ideas clearly. Indeed, there
are possible improvements in various different directions, which are currently work in progress.
For example, the numerical discretization we have used here is sub-optimal, particularly for
pattern-forming problems with spike-layers. One could investigate problems with much more 
complicated noise and/or aim to answer questions directly arising in larger-scale SPDE models 
from applications. (N1)-(N3) could also be applied to SPDEs not arising from forcing, but where 
coefficient functions are random variables. Furthermore, analyzing the stochastic fluctuations
numerically very close to bifurcation points requires new ideas to treat the full nonlinear 
stochastic system in this regime. Overall, the approach to focus 
directly on dynamical structures for SPDEs in combination with parameter continuation
could become another building block in uncertainty quantification.\medskip

\textbf{Acknowledgments:} I would like to thank the Austrian Academy of Sciences ({\"{O}AW})
for support via an APART fellowship. I also acknowledge support of the European Commission 
(EC/REA) via a Marie-Curie International Re-integration Grant. Furthermore, I would like to
acknowledge the helpful comments of two anonymous referees, which lead to improvements in the
presentation of the results.

\section{The Cubic-Quintic Allen-Cahn Equation}
\label{sec:AC}

Let $\cD\subset \R^d$ denote a bounded domain. As outlined in the introduction, 
our basic test problem will be based upon perturbing a version of the Allen-Cahn \cite{AllenCahn} 
equation. In general, the parabolic evolution Allen-Cahn PDE for $u=u(z,t)$ is given by 
\be
\label{eq:ac_evol}
\frac{\partial u}{\partial t}=\Delta u-\cF'(u;\mu),
\qquad (z,t)\in\cD\times [0,\I),\quad u:\cD\times [0,\I)\ra \R,
\ee
where $\Delta =\sum_{k=1}^d\frac{\partial^2}{\partial z_k^2}$ is the usual Laplacian, 
$\cF:\R\times \R\ra \R$ is the potential, $\mu\in\R$ is a parameter, $\cF'=:F$
indicates the derivative of $\cF$ with respect to the $u$-component 
and \eqref{eq:ac_evol} has to be augemented with suitable inital and boundary 
conditions. We remark that the general Allen-Cahn 
PDE \eqref{eq:ac_evol} is a well-studied system for various potentials given by 
polynomial nonlinearities. There has been a wide range of different analytical and 
numerical methods proposed to study Allen-Cahn-type PDEs, here we just mention a few 
examples \cite{AlfaroHilhorstMatano,AlmaBronsardGui,RabinowitzStredulinsky,Ward2}.
For the classical quartic potential leading to a cubic nonlinearity, the Allen-Cahn PDE 
is also known as the (real) Ginzburg-Landau PDE \cite{Gurtin,AransonKramer}, which is 
an amplitude equation or normal form for bifurcations 
from homogeneous states \cite{Mielke1,Schneider2}. Furthermore, the Ginzburg-Landau 
PDE is sometimes referred to as the Nagumo equation \cite{Nagumo,McKean,Popovic} developed 
as a simplified description for waves in nerve axons. Hence, the choice of deterministic 
Allen-Cahn PDE with a relatively general nonlinearity can definitely be viewed as a standard test 
problem for new methods. In this paper, we shall restrict to the case 
\be
d=2.
\ee
Numerically, the case $d=1$ can also be covered with the techniques developed here, 
while $d\geq 3$ will probably require more specialized numerical ideas. Basically we focus on 
$d=2$ here to have a sufficiently challenging numerical setup. The stationary problem 
associated to \eqref{eq:ac_evol}, for $u:\R^2\ra \R$, $u=u(z)$ with Dirichlet boundary 
conditions, is given by
\be
\label{eq:ac_det}
\left\{
\begin{array}{rcll}
0&=&\Delta u-F(u;\mu),\quad & \text{for $z\in\cD$,}\\
0&=&u,\quad & \text{for $z\in\partial \cD$.}\\
\end{array}
\right.
\ee
As a potential $\cF$ we are going to fix the following function
\be 
\label{eq:qc}
\cF(u;\mu):=2\mu u^2+u^4-\frac23 u^6 \qquad \Rightarrow \quad
F(u;\mu)=\cF'(u;\mu)=4(\mu u+u^3-u^5).
\ee
Since the nonlinearities in $F$ are of cubic and quintic types, one also 
refers to \eqref{eq:ac_evol} or \eqref{eq:ac_det} with the choice \eqref{eq:qc}
as the cubic-quintic Allen-Cahn (cqAC) equation. In fact, it is also very common 
to add a quintic term to the Ginzburg-Landau case, see for example 
\cite{DeisslerBrand,KapitulaSandstede,vanSaarloosHohenberg}, which again justifies
our choice of the qcAC equation as a basic test example.

\section{Numerical Continuation Basics}
\label{sec:numcont_basics}

It is natural to ask, which solutions the stationary cubic-quintic Allen-Cahn
equation \eqref{eq:ac_det} has under the variation of $\mu$. Numerical continuation
methods are a natural tool to address this problem. First, one uses a spatial 
discretization, {e.g.}, finite elements or finite differences to convert \eqref{eq:ac_det}
into a nonlinear algebraic system of equations
\be
\label{eq:ellipticPDE_cont1}
M(Z,\mu)=0,\qquad M:\R^{N_c}\times \R\ra \R^{N_c}, 
\ee
for some $Z\in\R^{N_c}$, where ${N_c}$ is usually very large for PDE problems. The idea of numerical 
continuation is to trace out a curve $\Gamma(s):=(\gamma(s),\mu(s))\in \R^{N_c}\times \R$ of solutions 
where $s$ parametrizes the the curve, for example by arclength. There are theoretically 
well-established predictor-corrector schemes available \cite{Govaerts,Seydel} 
to compute the curve $\Gamma(s)$. The main idea is that knowing a solution point $\Gamma(s_0)$ then one can 
compute $\Gamma(s_1)$ with $|s_0-s_1|=:\delta $, for $\delta>0$ sufficiently small, very efficiently. 
Indeed, first one uses a prediction step, for example by tangent prediction 
\be
\Gamma^{(0)}(s_1)=\Gamma(s_0)+\delta ~ \dot{\Gamma}(s_0),
\ee
where $\dot{\Gamma}(s_0)$ is the tangent vector to the curve $\Gamma(s)$ at $s_0$. Then ones uses 
a correction step, {i.e.}, one has to solve 
\be
\label{eq:ellipticPDE_cont2}
M(Z,\mu(s_1))=0,\qquad M:\R^{N_c}\times \R\ra \R^{N_c}, 
\ee
iteratively with starting initial guess $\Gamma^{(0)}(s_1)=(\gamma^{(0)}(s_1),\mu^{(0)}(s_1))$. If the
initial guess is sufficiently close to a true solution, which can be guaranteed under 
suitable non-degeneracy hypotheses, then an iterative solver converges very quickly to 
$\Gamma(s_1)=(\gamma(s_1),\mu(s_1))$. At degenerate bifurcation points, one may use branch switching 
techniques to track the bifurcating curves $\Gamma(s)$ beyond this point \cite{Kuznetsov}.\medskip

\section{Deterministic Continuation for the cqAC Equation}
\label{sec:det_cont}

\begin{figure}[htbp]
\centering
\psfrag{uu}{$\|u\|_2$}
\psfrag{p}{$\mu$}
\psfrag{x}{\scriptsize $x$}
\psfrag{y}{\scriptsize $y$}
\psfrag{u}{\scriptsize $u$}
\psfrag{a}{(a)}
\psfrag{b}{(b)}
\psfrag{c}{(c)}
\psfrag{d}{(d)}
\psfrag{bu}{\scriptsize (c)}
\psfrag{cu}{\scriptsize (d)}
\psfrag{au}{\scriptsize (b)}
\psfrag{R0}{$R_0$}
\psfrag{R1}{$R_1$}
\psfrag{R2}{$R_2$}
\psfrag{R3}{$R_3$}
\psfrag{R4}{$R_4$}
		\includegraphics[width=1\textwidth]{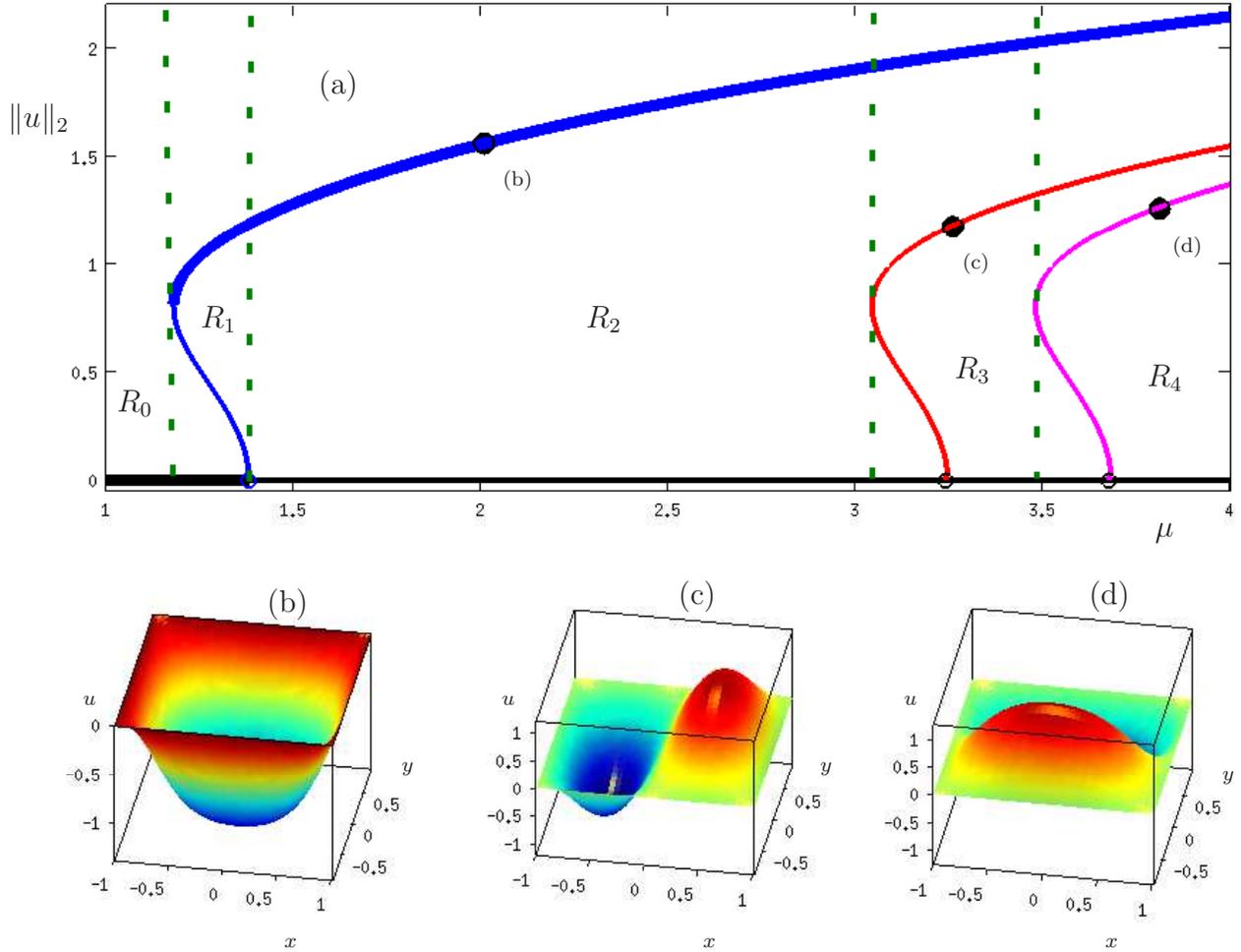}
\caption{\label{fig:1}(a) Bifurcation diagram of the cqAC equation \eqref{eq:ac_det}-\eqref{eq:qc} 
obtained from four continuation runs in \texttt{pde2path}; the vertical axis shows the 
$L^2(\cD)$-norm of the solution. The thick curves (black and blue)
are locally stable, {i.e.}, for the linearized problem the spectrum is contained in the left-half
of the complex plane. The three simple branch/bifurcation points \eqref{eq:bpoints} are marked with
circles. There are also three fold points \eqref{eq:fpoints} on the non-trivial bifurcation 
branches (blue/red/magenta). The five formally defined parameter regimes \eqref{eq:pregimes} 
are delineated by vertical dashed green lines. The three black dots indicate to which points 
the solutions in diagrams (b)-(d) are associated. (b)-(d) Solution plots for parameter values 
$\lambda=2.0075, 3.2579, 3.8093$ on the upper parts of the non-trivial branches.}	
\end{figure}

Figure \ref{fig:1} shows the results of numerical continuation for the cqAC equation 
\eqref{eq:ac_det}. The results were obtained from a standard example in the continuation
software \texttt{pde2path} \cite{UeckerWetzelRademacher} on a rectangle
\be
\cD=[-1,1]\times [-0.9,0.9].
\ee
The continuation was started with the homogeneous solution $u\equiv 0$, which 
exists for all $\mu\in\R$. Increasing $\mu$ from the starting value $\mu(s_0)=0$,
we focused on the first three bifurcation points. Those points are approximately obtained 
at 
\be
\label{eq:bpoints}
\mu^\txtb_1\approx 1.3798,\qquad \mu^\txtb_2\approx 3.2385,\qquad \mu^\txtb_3\approx 3.6784. 
\ee
We remark that it is known for the square domain $[-1,1]^2$ that the bifurcation points 
$\mu^\txtb_2$ and $\mu^\txtb_3$ coincide, which is one reason to consider the rectangle
$\cD$ to split up this numerical difficulty.

Up to $\mu^\txtb_1$, the homogeneous branch $\Gamma_0$ is locally asymptotically stable for the 
time-dependent cqAC equation \eqref{eq:ac_evol}. The linearized system along the homogeneous branch 
has at least one eigenvalue with positive real part for $\mu>\mu^\txtb_1$ and only eigenvalues with 
negative real part for $\mu<\mu^\txtb_1$. Branch switching at $\mu^\txtb_k$, $k\in\{1,2,3\}$, in 
\texttt{pde2path} leads to three subcritical bifurcation branches $\Gamma_k$, which each re-acquire 
one negative eigenvalue at further fold points located at 
\be
\label{eq:fpoints}
\mu^\txtf_1\approx 1.1794,\qquad \mu^\txtf_2\approx 3.0422,\qquad\mu^\txtf_3\approx 3.4779.
\ee
In particular, the branch $\Gamma_1$ restabilizes and the system is bistable for 
$\mu\in(\mu_1^\txtf,\mu_1^\txtb)$. We note the numerical results clearly indicate that 
$u\equiv 0$ is globally stable for $\mu<\mu_1^\txtf$, while the top part of the branch $\Gamma_1$
is the only attractor for $\mu>\mu_1^\txtb$.\medskip 

{\small \textbf{Remark:} At the bifurcation points, symmetries of the cqAc equation considered
here ($u\mapsto -u$ or $x_1\mapsto -x_1$ or $x_2\mapsto -x_2$) imply the existence of additional
bifurcation branches with symmetric solutions not shown in Figure \ref{fig:1}. In particular, 
throughout this manuscript we only select and analyze numerically one of the symmetric branches.
For more details on symmetry-breaking and using numerical continuation to find starting solutions
without homotopy continuation see, {e.g.},~\cite{KuehnEllipticCont}.}\medskip

We partition the parameter space into five 
disjoint intervals, which we are going to focus on throughout this manuscript
\be
\label{eq:pregimes}
R_0:=[0,\mu_1^\txtf),\quad 
R_1:=[\mu_1^\txtf,\mu_1^\txtb),\quad 
R_2:=[\mu_1^\txtb,\mu_2^\txtf),\quad 
R_3:=[\mu_2^\txtf,\mu_3^\txtf),\quad 
R_4:=[\mu_3^\txtf,4],\quad 
\ee
Basically, the regime $R_0$ is globally monostable, $R_1$ is bistable, $R_3$ has one stable and 
one simple homogeneous unstable branch, while $R_3$ and $R_4$ consider the interaction between 
one and two non-trivial unstable branches $\Gamma_2$ and $\Gamma_3$, respectively, with the stable branch 
$\Gamma_1$ and the homogeneous branch $u\equiv 0$. 

\section{The Stochastic Allen-Cahn Equation}
\label{sec:AC_SPDE_intro}

We return to the general Allen-Cahn equation \eqref{eq:ac_evol}. One natural approach to include noise
is to replace the PDE, on a very formal level, by a stochastic partial differential equation (SPDE) 
given by
\be
\label{eq:SPDE_AC1}
\frac{\partial u}{\partial t}=\Delta u-F(u;\mu)+G(u)\eta,
\ee
where $\eta=\eta(x,t)$ is a space-time-dependent stochastic process and $G$ is a suitable mapping. 
The Allen-Cahn SPDE has gained - unsurpringly considering the discussion in 
Section \ref{sec:AC} - considerable attention in various analytical studies, 
see {e.g.}~\cite{FarisJona-Lasinio,BerglundGentz10,BloemkerMohammed1,HairerRyserWeber,
KohnOttoReznikoffVanden-Eijnden,Weber}. 
Furthermore, numerical methods for SPDEs have frequently been tested on 2D Allen-Cahn-type SPDEs and 
related equations \cite{CastroLythe,KovacsLarssonLindgren,LordThuemmler,RyserNigamTupper,Shardlow}. 
Hence, \eqref{eq:SPDE_AC1} can also be viewed as a standard test problem for the SPDE case.\medskip

To make the formal SPDE \eqref{eq:SPDE_AC1} precise, we shall adopt the framework of mild solutions considerd 
in \cite{DaPratoZabczyk}. Consider a probability space $(\Omega,\cG,\P)$ and a Hilbert space $\cH$; for 
practical purposes we may always think of $\cH=L^2(\cD)$ or $\cH=H_0^1(\cD)$ here. Denote 
by $\langle\cdot,\cdot\rangle_{\cH}$ the inner product in $\cH$ and by 
$\|\cdot\|_\cH=\sqrt{\langle\cdot,\cdot\rangle_{\cH}}$ the associated norm. Let $Q:\cH\ra \cH$ 
be a symmetric non-negative 
linear operator and suppose $Q$ is of trace-class $\textnormal{Tr}(Q)<+\I$. Then there exists a 
complete orthonormal system $\{e_k\}_{k=1}^\I$ for the Hilbert space $\cH$ such that
\benn
Qe_k=\lambda_k e_k,\qquad \text{for $k\in\N$,}
\eenn
where $\{\lambda_k\}_{k=1}^\I$ is a non-negative bounded summable of eigenvalues. Then one may 
define a $Q$-Wiener process by
\be
\label{eq:series}
W(t):=\sum_{k=1}^{\I}\sqrt{\lambda_k}\beta_k(t)e_k
\ee
where $\{\beta_k(t)\}_{k=1}^\I$ are independent one-dimensional Brownian motions. The series 
\eqref{eq:series} converges in $L^2(\Omega,\cG,\P;\cH)$; see also \cite[p.86-89]{DaPratoZabczyk}. 
One finds that $\E[W(t)]=0$ and for all $f_1,f_2\in \cH$ we have $\E[\langle f_1,W(t)\rangle 
\langle f_2,W(s)\rangle ]=\min(t,s)\langle Qf_1,f_2\rangle$ so that $Q$ can be viewed as the covariance 
operator. Furthermore, define $\cH_0:=Q^{1/2}\cH$ and let $G:\cH\ra L_2(\cH_0,\cH)$,
where $L_2(\cH_0,\cH)$ denotes the space of Hilbert-Schmidt operators. In this case one may consider 
the stochastic Allen-Cahn equation for $u(t)=u(\cdot,t)\in \cH$ as an $\cH$-valued stochastic 
evolution equation
\be
\label{eq:AC_SPDE}
\txtd u=\left[\Delta u-F(u;\mu)\right]\txtd t+G(u)~\txtd W,\qquad u_0:=u(0)\in \cH,\\
\ee
which can interpreted in the usual precise integral form \cite{DaPratoZabczyk} by stochastic integrals;
note that in \eqref{eq:AC_SPDE}, we understand $F$ in the Nemytskii operator 
sense as discussed in \cite[Sec.7.2]{DaPratoZabczyk} and mild solutions exist under suitable 
local Lipschitz conditions on $F$. Furthermore, we shall always assume that 
the initial condition is deterministic and lies in a suitable subspace of 
$\cH$, e.g.~$u(0)\in C(\cD)\subset L^2(\cD)$ as discussed in \cite{DaPratoZabczyk}.\medskip

The solution theory for general SPDEs is quite involved and is still 
being developed for certain cases. However,
we are going to work with a linearized version of \eqref{eq:AC_SPDE} for which a suitable theory
of mild solutions exists; see Section \ref{sec:abs_Lyapunov}.

\section{The Abstract Lyapunov Equation}
\label{sec:abs_Lyapunov}

We have seen in Section \ref{sec:det_cont} that the deterministic cqAC PDE has several locally
asymptotically stable stationary solutions. It is then a natural question to ask, what influence
different forms of the noise term $G(u)~\txtd W$ have on the fluctuations around these stationary
solutions. Here we are interested in fluctuations on a subexponential time-scale, {i.e.}, before
the regime of large deviation theory \cite{FreidlinWentzell} is reached. This implies that we 
restrict to either
\begin{itemize}
 \item[(R1)] sufficiently small norm of $G(u)$, or
 \item[(R2)] sufficiently small eigenvalues $\lambda_k$, or
 \item[(R3)] sufficiently small time final $T$ such that no large deviations have occurred. 
\end{itemize}
In this regard, it is important to observe that (R1)-(R3) are related. For example, may absorb 
a constant positive multiplicative prefactor $\zeta_1>0$ in $\zeta_1~G(u)~\txtd W$ either
into $G$ by defnining $G_{\textnormal{new}}:=\zeta_1G$, or into the eigenvalues $\lambda_k$ 
associated to $Q$ by defining a new operator $Q_{\textnormal{new}}$. Furthermore, one may apply
a change of time scale in \eqref{eq:AC_SPDE} via $t_{\textnormal{new}}:=\zeta_2 t$ for a new
time $t_{\textnormal{new}}$ and a positive scale factor $\zeta_2>0$ to relate the condition (R3)
to (R1)-(R2).\medskip

In one (or more) of the cases (R1)-(R3), it natural
to consider a linearized form of \eqref{eq:AC_SPDE} around a deterministic stationary
solution $u^*=u^*(z)$ for a parameter value $\mu=\mu^*$. The linear SPDE is given by
\be
\label{eq:SPDE_lin}
\txtd U=A U ~\txtd t+B~\txtd W,\qquad U_0:=U(0)\in\cH,\\
\ee
where $U(t)=U(\cdot,t)\in\cH$, the maps $A,B:H\ra H$ are linear operators given by
\be
A=\Delta - \cF''(u^*;\mu^*),\qquad B=G(u^*)
\ee
and $\cF''$ denotes the second derivative of $\cF$ with respect to the $u$-component, {i.e.}, we have
\be
\cF''(u;\mu)=4(\mu+3(u^*)^2-5(u^*)^4)\textnormal{Id}.
\ee
where $\textnormal{Id}:\cH\ra \cH$ just denotes the identity operator. Note carefully that the
spectrum of $A$, denoted here by $\textnormal{spec}(A)\subset \C$, indicates whether $u^*$ is a locally 
stable solution of the deterministic cqAC PDE \eqref{eq:ac_evol}. If 
$\textnormal{spec}(A)$ is properly contained in the left-half complex plane, then we have 
local asymptotic stability.
This is the main case we are going to focus on here; recall that the locally stable stationary 
solutions of the cqAC PDE have been computed solving \eqref{eq:ac_det} by numerical continuation 
in Section \ref{sec:det_cont} and are shown as the thick curves in Figure \ref{fig:1}.\medskip

The linear SPDE \eqref{eq:SPDE_lin} is the direct analog of the stochastic ordinary 
differential equation for an Ornstein-Uhlenbeck (OU) process. As a solution concept, we consider 
mild solutions \cite[{Ch.7}]{DaPratoZabczyk} given by
\be
\label{eq:mild}
U(t)=S(t)U_0+\int_0^tS(t-s)B~\txtd W(s).
\ee
where $S(t)=\txte^{tA}$ is the semigroup generated by $A$. It is helpful to denote the stochastic 
integral in \eqref{eq:mild} as
\benn
W_{A,B}(t):=\int_0^tS(t-s)B~\txtd W(s).
\eenn
$W_{A,B}(t)$ is also called the stochastic convolution, which has the associated covariance operator 
\cite[{Thm. 5.2}]{DaPratoZabczyk}
\be
\label{eq:CovConv}
V(t):=\textnormal{Cov}( W_{A,B}(t))=\int_0^t S(r)BQB^*S^*(r)~\txtd r
\ee
where $B^*$ denotes the adjoint operator of $B$ and $S^*(r)$ denotes the adjoint semigroup of 
$S(r)$. Computing $V(t)$ provides the main relevant information on the local fluctuations 
of the OU-process. Since the OU process 
differs only in the quadratic and higher-order terms from the local dynamics around the 
stationary solution $u^*$ of the cqAC SPDE, one may show that on subexponential time scales, the
operator $V(t)$ provides the leading-order fluctuations around $u^*$; see also Appendix \ref{ap:2}. Using 
\eqref{eq:series}, one can write the OU process solution also as
\be
\label{eq:formal_series_sol}
U(t)=S(t)U_0+\sum_{k=1}^\I\sqrt{\lambda_k}\int_0^t S(t-s)Bq_k~\txtd \beta_k(s),
\ee
where the detailed definition of the integral is given in \cite[{Sec. 2.2}]{DaPrato1}. 
Using a standard result \cite[{Prop. 2.2}]{DaPrato1}, it is not difficult to see that 
the series in \eqref{eq:formal_series_sol} is convergent in $L^2(\Omega,\cG,\P;\cH)$. 
Furthermore, $U(t)$ is a Gaussian random variable with mean $S(t)U_0$ and covariance operator 
$V(t)$. Indeed, if we are in the case where $u^*$ is locally stable for the deterministic 
cqAC PDE then we have
\benn
\|S(t)\|_\cH\leq K \txte^{\rho t}
\eenn
for some constant $K>0$ and $\rho<0$, where 
$\textnormal{spec}(A)\subset \{\gamma \in\C:\text{Re}(\gamma)<\rho\}$. Therefore, 
$\|S(t)U_0\|_\cH\ra 0$ as $t\ra +\I$ exponentially fast, so we may formally 
neglect the first term $S(t)U_0$. In this case, we simply have $U(t)=W_{A,B}(t)$ and 
there exists a unique invariant Gaussian measure with mean zero and covariance operator 
$V_\I:=\lim_{t\ra +\I} V(t)$ \cite[{Thm. 2.34}]{DaPrato1}; again the convergence to $V_\I$ is
expected to be exponentially fast so we may negelect the initial fast transients and just 
focus on $V_\I$. $V_\I$ is a linear continuous operator that satisfies 
$\langle V_\I h,h\rangle_{\cH}\geq 0$ for all $h\in \cH$ and $V_\I=V_\I^*$. Furthermore, it is 
crucial to note that $V_\I$ satisfies a Lyapunov equation \cite[{Lem. 2.45}]{DaPrato1} given by
\be
\label{eq:Lyapunov}
\langle V_\I g,A^*h\rangle_{\cH} +\langle A^*g,V_\I h\rangle_{\cH}=-\langle BQB^* g,h\rangle_{\cH},
\ee
for any $g,h\in \cH$. Observe that \eqref{eq:Lyapunov} has the same algebraic structure of 
the classical Lyapunov equation associated with linear stochastic ordinary differential 
equations (SODEs) \cite[Ch.5]{BerglundGentz}. The result \eqref{eq:Lyapunov} can also be extended to apply in 
Banach spaces \cite[{Sec.4}]{GoldysvanNeerven}.

\section{A Finite-Difference Galerkin Scheme}
\label{sec:fd}

One possibility to compute the operator $V_\I$ is to start directly from the abstract 
Lyapunov equation \eqref{eq:Lyapunov} and derive suitable numerical methods for it. 
However, this studies \eqref{eq:Lyapunov} in isolation from the other numerical 
algorithms one may develop for the deterministic and stochastic cqAC equation. Here we 
aim to determine a finite-dimensional approximation to $V_\I$ by starting 
from the SPDE \eqref{eq:AC_SPDE}.\medskip

As a first step, we use a spatial finite-difference discretization for the SPDE \eqref{eq:AC_SPDE}
based upon the one-dimensional scheme presented in \cite{SauerStannat1}, in combination with
a spectral approximation of the noise \cite{Shardlow}. The reason for this combination will be 
explained below, after we have introduced the scheme. Let 
$\{(x_m,y_n)\}\subset \R^2$ for $m\in\{0,1,2,\ldots,M\}$ and $n\in\{0,1,2,\ldots,N\}$ denote 
the points in a regular spatial mesh of edge lengths $h_x,h_y$ for the rectangle 
$\Omega=[-L_x,L_x]\times [-L_y,L_y]$, {i.e.}, we take
\be
\begin{array}{lll}
x_m=-L_x+m h_x,\quad & m\in\{0,1,2,\ldots,M\},\qquad & h_x:=2L_x/M,\\ 
y_n=-L_y+n h_y,\quad & n\in\{0,1,2,\ldots,N\},\qquad & h_y:=2L_y/N,\\
\end{array}
\ee
which yields $(M+1)(N+1)$ points and $(M-1)(N-1)$ interior points. Denote the approximation 
of $u$ on the mesh by $u^{m,n}\approx u(x_m,y_n)$. We approximate the Laplacian on the interior points, 
assuming Dirichlet boundary conditions, by the usual five-point stencil given by
\benn
(\Delta u)(x_m,y_n)\approx\frac{u^{m+1,n}+u^{m-1,n}-2u^{m,n}}{(h_x)^2}
+\frac{u^{m,n+1}+u^{m,n-1}-2u^{m,n}}{(h_y)^2}=:(\Delta u)^{m,n}.
\eenn 
The nonlinear terms are approximated by evaluation at the mesh points
\be
F^{m,n}:=F(u(x_m,y_n);\mu),\qquad G^{m,n}:=G(u(x_m,y_n)), 
\ee
where the $\mu$-dependence of the deterministic nonlinear part will be omitted in the
notation for simplicity. Next, denote by 
\be
\label{eq:indicator}
{\bf 1}^{m,n}(z):=\frac{1}{h_x h_y}{\bf 1}_{[x_m,x_{m+1})\times [y_n,y_{n+1})}(z)
\ee
the $L^2(\cD)$-normalized indicator function of the rectangle $[x_m,x_{m+1})\times [y_n,y_{n+1})$. 
To approximate the covariance operator $Q$ consider, the expansion      
\be
\label{eq:series_approx}
W(t)=\sum_{k=1}^{\I}\sqrt{\lambda_k}\beta_k(t)e_k\approx 
\sum_{k=1}^{K}\sqrt{\lambda_k}\beta_k(t)e_k,
\ee
where we may select $K$ to control the approximation accuracy of the noise term. Then
we consider the inner product of the cqAC SPDE with the indicator function 
\eqref{eq:indicator} to obtain
\be
\txtd \langle u,{\bf 1}^{m,n} \rangle_{L^2} =\left[\langle \Delta u,{\bf 1}^{m,n}\rangle_{L^2}
-\langle F(u;\mu),{\bf 1}^{m,n}\rangle_{L^2}\right]\txtd t+ 
\sum_{k=1}^{\I} \sqrt{\lambda_k} \langle G(u)e_k,{\bf 1}^{m,n}\rangle_{L^2}~ \txtd\beta_k(t).
\ee
Now one may consider each of the different terms for $\cH=L^2(\cD)$ and note that
\be
\langle u,{\bf 1}^{m,n} \rangle_{L^2}=\int_\Omega u(z) ~{\bf 1}^{m,n}(z)~\txtd z
\approx u^{m,n}\frac{h_x h_y}{h_x h_y}=u^{m,n}. 
\ee
Similarly, we have
\be
\langle \Delta u,{\bf 1}^{m,n}\rangle_{L^2} \approx (\Delta u)^{m,n},\qquad 
\langle F(u;\mu),{\bf 1}^{m,n}\rangle_{L^2} \approx F^{m,n}.
\ee
For the noise term, we use a truncation level $K$ as in \eqref{eq:series_approx}, and 
obtain
\be
\sum_{k=1}^{\I} \sqrt{\lambda_k} \langle G(u)e_k,{\bf 1}^{m,n}\rangle_{L^2}~ \txtd\beta_k(t)
\approx \sum_{k=1}^{K} \sqrt{\lambda_k} G^{m,n} e_k^{m,n} \txtd\beta_k(t),
\ee
where $e_k^{m,n}:=\langle e_k,{\bf 1}^{m,n}\rangle$ are scalars which can be pre-computed
in an offline step. We use the notation $p_j:=u^{m,n}$, where $j$ is an index 
labelling the interior vertices of the mesh and which ranges between $1$ and $(M-1)(N-1)$, consider
the vector $p:=\{p_j\}_{j=1}^{(M-1)(N-1)}\in\R^{(M-1)(N-1)}=:\R^J$ and define
\be
\vartheta_j(p;\mu):=(\Delta u)^{m,n}-F^{m,n},\qquad \sigma_{k,j}(v):= \sqrt{\lambda_k}G^{m,n} e_k^{m,n}.
\ee
so that $\vartheta:\R^J\times \R\ra\R^J$ is just a vector-valued map and 
$\sigma:\R^{J}\ra \R^{K\times J}$ is a matrix-valued map. Hence, we may approximate the 
SPDE with a $J$-dimensional system of SODEs given by
\be
\label{eq:SODE}
\txtd p=\vartheta(p;\mu)~ \txtd t+\sigma(p)~\txtd \beta,
\ee
where $\beta=\beta(t)=(\beta_1(t),\beta_2(t),\ldots,\beta_K(t))^\top \in \R^K$ is a 
vector of $K$ independent identically distributed (iid) Brownian motions. Note that 
we have now two different ways to
control the approximation level. The spatial mesh can be refined by increasing $M$ 
and/or $N$ while the spectral approximation of the noise can be improved by 
increasing $K$. Of course, many other spatial discretization schemes would lead 
to a system of SODEs. The next step starts from the level of the SODEs \eqref{eq:SODE}
to approximate the covariance operator of the problem.

\section{Lyapunov Equations for SODEs}
\label{sec:LyaSODEs}

We start from the SODE \eqref{eq:SODE}, where we are interested primarily in the case
where the noise level is suitably bounded, {i.e.}, $|\sigma(p)|< \sigma^*$ for all $p$ in
some compact set of phase space and $\sigma^*$ is a fixed constant. This is quite a natural 
viewpoint as very large noise terms may indicate that the model might have not been derived 
to suitable physical modelling accuracy in the first place. Here we briefly recall the 
idea to approximate certain covariance operators associated to \eqref{eq:SODE} along 
bifurcation curves developed in \cite{KuehnSDEcont1}.\medskip 

Suppose that the associated ODE for $\sigma\equiv 0$, given by 
\be
\label{eq:ODE}
\dot{p}=\frac{\txtd p}{\txtd t}=\vartheta(p;\mu),
\ee
has a hyperbolic stable steady state $p^*=p^*(\mu)$ for a given range of parameter values. 
For the cqAC SPDE case we consider here, we may think of the discretized versions of the 
stable parts of the branches $\Gamma_0$ and $\Gamma_1$ computed in Section \ref{sec:det_cont}.
For example, for the homogeneous branch $\Gamma_0$ in the parameter range $\mu \in R_0\cup R_1$
we have the stable solution $(u^*)^{m,n}=(p^*)^j= 0$ for all $j$. If we suppose that 
$\sigma(p^*)\neq 0$ and Taylor expand \eqref{eq:SODE} to lowest order at $p^*$ we obtain
\be
\label{eq:SDE2}
\txtd P=\cA P ~\txtd t  + \cB~ \txtd \beta, \qquad P=P(t)\in\R^J,
\ee
where $\cA=\cA(p^*;\mu)=(D_p\vartheta)(v^*;\mu)\in\R^{J\times J}$ is the usual Jacobian matrix
and $\cB=\sigma(p^*)\in \R^{J\times K}$ is the lowest order term of the noise. Equation 
\eqref{eq:SDE2} is a $J$-dimensional Ornstein-Uhlenbeck (OU) process \cite{Gardiner} with 
covariance matrix 
\begin{equation*}
\cV(t):=\text{Cov}(P(t))=\int_0^t \cU(t,s)\cB\cB^\top \cU(t,s)^\top~ \txtd s
\end{equation*}
where $\cU(t,s)$ is the fundamental solution of $\dot{\cU}=\cA\cU$. Differentiation shows 
that $\cV(t)$ satisfies the ODE 
\be
\label{eq:Cov_ODE}
\frac{\txtd\cV}{\txtd t}=\cA \cV+\cV\cA^\top+\cB\cB^\top. 
\ee
Since $p^*$ is a hyperbolic stable equilibrium point, it follows \cite{Bellman} that the 
eigenvalues of the linear operator $\cL(\cV):=\cA \cV+\cV\cA^\top$ are $\{2\lambda_j\}_{j=1}^J$, 
where $\lambda_j$ are the eigenvalues of $\cA$. Therefore, the covariance ODE \eqref{eq:Cov_ODE} 
has a globally stable equilibrium solution defined by 
\be
\cV_\I:=\lim_{t\ra +\I} \cV(t), 
\ee
which is obtained by solving
\be
\label{eq:Lya1}
0=\cA \cV_\I+\cV_\I\cA^\top+\cB\cB^\top.
\ee
Note carefully that \eqref{eq:Lya1} is a Lyapunov equation and $\cV_\I$ is a finite-dimensional
approximation to the covariance operator $V_\I$ discussed in Section \ref{sec:abs_Lyapunov}. 
We denote the solution of \eqref{eq:Lya1} sometimes also as $\cV_\I(p^*(\mu);\mu)$ to stress 
the dependence upon the parameter $\mu$ and the associated state $p^*(\mu)$.\medskip 

It is standard \cite{BerglundGentz} to define a covariance neighbourhood
\be
\label{eq:Bh_big}
\cC(r):=\left\{p\in\R^n:(p-p^*)^\top\cV_\I^{-1}(p-p^*)\leq r^2\right\}
\ee
where $r$ is a parameter that can be interpreted as a probabilistic confidence level. 
Essentially, $\cC(r)$ yields the information about local fluctuations and important 
noise directions near $p^*$. Hence, we would like to solve for $\cV_\I$ along the entire branch 
of steady states obtained via continuation $\Gamma:=\{(p^*(\mu),\mu)\}\subset \R^J\times \R$. 
The algebraic equation \eqref{eq:Lya1} is a uniquely solvable Lyapunov matrix equation \cite{GajicQureshi}. 
If the dimension $J$ of the problem is moderate, direct numerical methods for \eqref{eq:Lya1} 
work quite well. A well-known direct solution algorithm is the Bartels-Stewart algorithm
\cite{BartelsStewart,GolubNashVanLoan}.\medskip 

However, we remark that for our case new aspects arise since we want to solve 
\eqref{eq:Lya1} along an entire steady state branch $\Gamma$. Indeed, most numerical 
continuation algorithms require an approximation of the $J\times (J+1)$ Jacobian matrix 
$(D_{(p,\mu)}\vartheta)(p^*(\mu_1),\mu_1)$ to compute a point $(p^*(\mu_2),\mu_2)\in\Gamma$ 
starting from $(p^*(\mu_1),\mu_1)\in \Gamma$. Therefore, the matrix 
$\cA=(D_p\vartheta)(p^*(\mu_1);\mu_1)$ is available at each continuation step. Furthermore, 
assembling the matrix $\cB$ at a given point $(p^*(\mu_1),\mu_1)$ is relatively cheap numerically
if $K$ is of moderate size. Solving \eqref{eq:Lya1} at $(p^*(\mu_1),\mu_1)$ gives a matrix 
$\cV_\I(p^*(\mu_1);\mu_1)$. If $|\mu_1-\mu_2|$ is small then $\cV_\I(p^*(\mu_1);\mu_1)$, or a 
tangent prediction starting from $\cV_\I(p^*(\mu_1);\mu_1)$, is already an excellent initial guess 
to find $\cV_\I(p^*(\mu_2);\mu_2)$. Hence, except for the first point on the equilibrium curve, we 
always have an initial guess available for an iterative method of solving the Lyapunov equation. 
This was exploited in \cite{KuehnSDEcont1} via a Gauss-Seidel scheme by re-writing the 
Lyapunov equation as a linear system. 

\section{Continuation of the Lyapunov Equation for cqAC}
\label{sec:contSPDE_Lya}

In this paper, we also re-write the Lyapunov equation \eqref{eq:Lya1} as a large (sparse) linear system 
using the Kronecker product $\otimes$ in the usual way
\be
\label{eq:Lya_linsyst}
\left[\text{Id} \otimes \cA+ \cA\otimes \text{Id}\right]\text{vec}(\cV_\I) =-\text{vec}(\cB\cB^\top),
\ee
where $\text{Id}$ is the identity matrix of size $J\times J$ and $\text{vec}(\cdot)$ converts a matrix 
$\cM$ to a vector by stacking it column-by-column, {e.g.}, 
$\text{vec}(\cM)=(\cM_{11}, \cM_{21},\ldots, \cM_{12}, \cM_{22},\ldots)^\top$. Here we are going to
use standard iterative solvers for sparse linear systems to continue the Lyapunov equation as implemented
in \texttt{MatLab}. In particular, we use \texttt{bicgstab} (a stabilized biconjugate gradient method, our 
default choice here), \texttt{gmres} (a generalized minimum residual method) and \texttt{qmr} (a quasi-minimal
residual method). We use as an initial guess the previous solution of the Lyapunov equation on the continuation 
branch; if we have no solution available, the zero vector is used as an initial guess. \medskip

To obtain \eqref{eq:Lya_linsyst}, we also have to make a choice of the mesh parameters described in
Section \ref{sec:fd}. Note that we work on a regular grid on the rectanglular domain 
\be
\cD=[-1,1]\times[-0.9,0.9]=[-L_x,L_x]\times[-L_y,L_y]
\ee
if and only if $h_x=2L_x/M\stackrel{!}{=}2L_y/N=h_y$, which implies upon fixing $M$ that
$N=ML_y/L_x=\frac{9}{10}M$. Hence, to obtain the case $h_x=h_y$, we must require that 
$\frac{9}{10}M$ is a natural number. One good choice is $M=50$ and so $N=45$, which will
be our standard choice here. Recall again that here we do not aim for optimal numerical 
performance of the algorithm. 
Indeed, it is important to note that the method we use here to continue the
Lyapunov equation is numerically sub-optimal in several different ways:

\begin{enumerate}
 \item[(S1)] There are more efficient methods for iterative solution of matrix Lyapunov equations 
 available that take direct advantage of the special matrix-equation structure, {e.g.}~ADI-type methods
 would be an option \cite{BennerLiPenzl,GajicQureshi}.
 \item[(S2)] Instead of using the previous solution as an initial guess, a full numerical continuation
 scheme may also employ other predictors, {e.g.}~a tangent predictor \cite{Govaerts}.
 \item[(S3)] Using a regular mesh is likely to be sub-optimal in many cases. One could use a suitable 
 finite-element method spatial discretization together with adaptive meshing for the spatial discretization, 
 which is going to improve the performance in many situations.
\end{enumerate}

If, despite using a sub-optimal numerical algorithm as explained in (S1)-(S3), we can still
compute the Lyapunov equation efficiently on a standard desktop hardware platform\footnote{Basic details 
of the desktop computer setup used: Intel Core i5-4430 CPU @ 3.00GHz processor (quadcore), Kingston 16GB system
memory, WDC WD10EZRX-00L 1TB harddrive, Ubuntu 12.04.4 LTS operating system, MatLab R2013a computational 
platform.} with our basic continuation setup, then we have achieved a valid \emph{proof-of-concept}. In fact,
this leads to the conjecture that the computation time can be reduced potentially by one or more orders
of magnitude if all available ideas in (S1)-(S3) are explored.\medskip

We use the domain $\cD$ and set $\cH:=H^1_0(\cD)$ as our Hilbert 
space. Furthermore, as a basic test noise we consider a decaying noise term with eigenvalues 
$\lambda_k$ for $Q$ given by
\be
\label{eq:noise_num}
\lambda_k=\tilde{\sigma} \exp\left(-\frac{1}{10} \varphi_k\right),\qquad k\in\{1,2,\ldots,K\},
\ee
where $\tilde{\sigma}>0$ is a parameter, which can be used to tune the noise level and 
$\varphi=(\varphi_1,\ldots,\varphi_K)^\top$ is a vector with increasing non-negative entries to be 
specified. As a basis $e_k$, we fix
\be
\label{eq:efuncs}
e_k(x,y)=\sin\left(\frac{\pi k_1}{2 L_x}(x+L_x)\right)\sin\left(\frac{\pi k_2}{2 L_y}(y+L_y)\right),
\qquad k=k_1+k_2,\quad k_1,k_2\in\N.
\ee
As default parameters for the continuation of the Lyapunov equations we set $\tilde{\sigma}=5$ and consider
an additive noise as a first test case so that $G(u)\equiv 1$. Regarding the linear system solvers, we 
fix the maximum iteration number for the iterative linear system solvers to $\texttt{maxit}=200$ and
the tolerance to $\texttt{tollin}=10^{-4}$.

\subsection{The Lower Branch}
\label{ssec:lower_res}

For the first numerical test case, we focus on the homogeneous stationary solution branch $u^*\equiv 0$, 
which is locally asymptotically stable for $\mu\in[0,\mu_1^b)$ with $\mu_1^b\approx 1.3798$ for the deterministic 
system. Three different noise terms were considered with eigenvalues $\lambda_k$ given by 
\eqref{eq:noise_num} associated to three different truncation levels $K=2,4,8$, more precisely, 
we have taken $\varphi_k=k$ for $k\in\{1,2,\ldots,K\}$. \medskip

\begin{figure}[htbp]
\centering
\psfrag{t}{$T$}
\psfrag{iter}{$\cI$}
\psfrag{mu}{$\mu$}
\psfrag{a}{(a)}
\psfrag{b}{(b)}
		\includegraphics[width=1\textwidth]{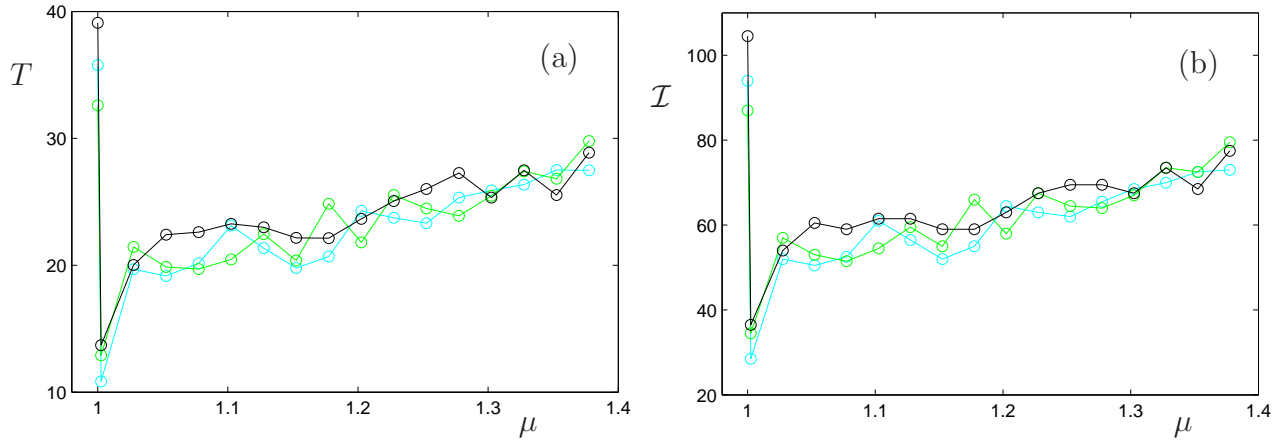}
\caption{\label{fig:2}The computation time $T$ and number of iterations $\cI$ for the Lyapunov operator
computation for the lower homogeneous branch with $\mu\in[0,\mu_1^\txtb)$ are shown. Each circle represents 
a point on the continuation curve,
where \eqref{eq:Lya1} has been solved using the \texttt{bicgstab}-method. The linear interpolation between
the points is for visualization purposes. (a) Parameter $\mu$ plotted against computation time
$T$ (in seconds). The cyan ($K=2$), green ($K=4$) and black ($K=8$) are different truncation
levels for the noise. (b) Number of iterations $\cI$ required for \texttt{bicgstab} to reach a tolerance of
less than $10^{-4}$. The color coding is as for (a).}	
\end{figure}

\begin{figure}[htbp]
\centering
\psfrag{lam}{$\lambda_k$}
\psfrag{k}{$k$}
\psfrag{un}{$\|\cdot\|$}
\psfrag{us+V}{$\|\cV_\I\|$}
\psfrag{us-V}{$-\|\cV_\I\|$}
\psfrag{us}{$\|u\|$}
\psfrag{mu}{$\mu$}
\psfrag{a}{(a)}
\psfrag{b}{(b)}
\psfrag{c}{(c)}
		\includegraphics[width=1\textwidth]{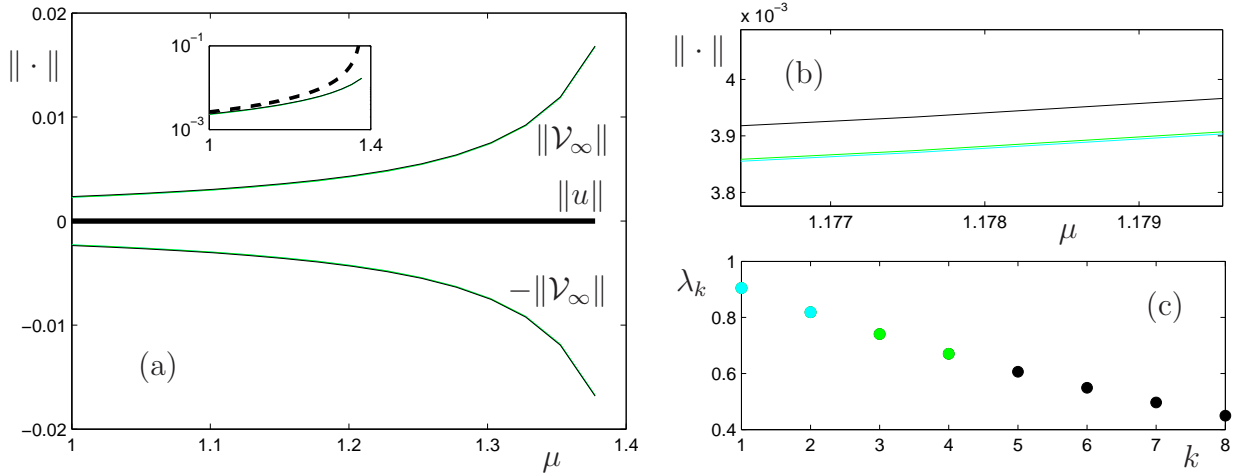}
\caption{\label{fig:3}(a) All norms in this figure are maximum norms $\|\cdot\|=\|\cdot\|_\I$. 
The deterministic homogeneous steady state solution $u^*\equiv 0$ 
(thick black curve) is shown, which also corresponds to the zero steady state $p^*\equiv 0$ of the discretized 
system. Furthermore, the norm of the covarince matrix $\cV_\I$ is plotted to indicate the size of the
fluctuations as measured by the linearized system; there are actually three different neighbourhoods 
shown in cyan ($K=2$), green ($K=4$) and black ($K=8$) corresponding to different truncation
levels for the noise. The inset in (a) shows a semi-logarithmic plot of the upper curves $\|\cV_\I\|$ (thin 
curves) in comparison to the scaling law $\cO(1/(\mu_1^b-\mu))$ (dashed black thick curve). (b) Zoom of 
part (a) near the upper part of the covariance neighbourhood to show that the curves given by 
$\|\cV_\I\|_\I$ for $K=2,4,8$ are very close. (c) Eigenvalues $\lambda_k$ used for the noise, {i.e.}, for 
$K=2$ there are two eigenvalues 
(in cyan), for $K=4$ there are fo$M=50$ and so $N=45$ur eigenvalues (the cyan ones and the green ones), while for $K=8$ also
the black ones are taken into account.}
\end{figure}

The results of the computations are shown in Figures \ref{fig:2}-\ref{fig:3}. Figure \ref{fig:2} considers 
the basic computational performance for the large sparse linear systems solution to find $\cV_\I$. The 
computational time $T$ and the number of iterations $\cI$ are very large in the first step, as we do
not have a good initial guess available and just start with the zero vector. Afterwards, each continuation
point takes approximately between 20 and 30 seconds on a standard desktop computer so the entire branch
of 17 points is processed within a few minutes. We conjecture that taking (S1)-(S3) into account, the
computation time can be reduced to a few seconds for the entire branch on the same hardware. In any case,
even for the sub-optimal computational framework, we can quickly analyze the local influence of noise and
quantify the uncertainty of the local system dynamics for different parameter values. 

\begin{figure}[htbp]
\centering
\psfrag{1norm}{$\|\cdot\|_1$}
\psfrag{2norm}{$\|\cdot\|_2$}
\psfrag{infnorm}{$\|\cdot\|_\I$}
\psfrag{mu}{$\mu$}
\psfrag{a}{(a)}
\psfrag{b}{(b)}
\psfrag{c}{(c)}
		\includegraphics[width=1\textwidth]{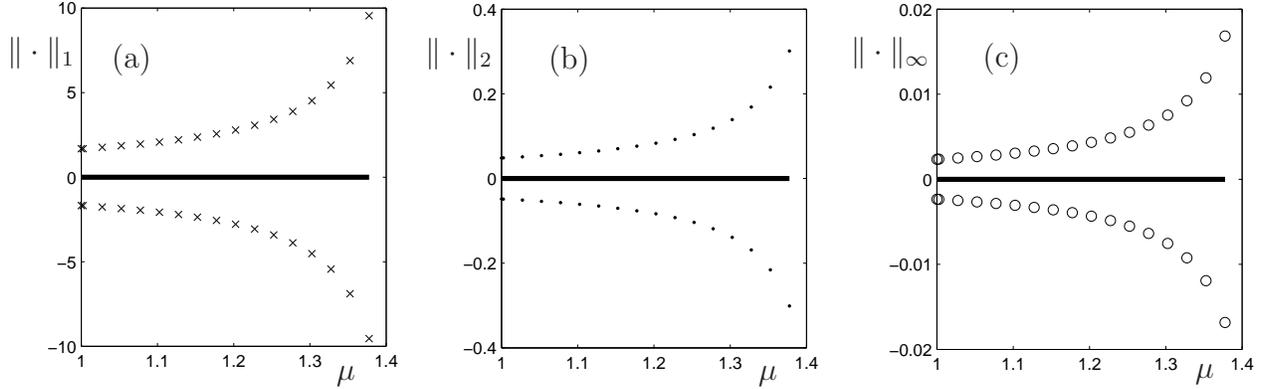}
\caption{\label{fig:4}Comparison of different norms for the vector of variances $\text{diag}(\cV_\I)$, which is 
obtained as the diagonal of $\cV_\I$. As usual, the horizontal axis in (a)-(c) shows the 
parameter values for $\mu$, where the homogeneous zero branch $u^*\equiv 0$ (thick black) is locally 
asymptotically stable. The figures (a)-(c), respectively crosses, dots and circles, correspond to 
$\|\text{diag}(\cV_\I)\|_1$, $\|\text{diag}(\cV_\I)\|_2$ and $\|\text{diag}(\cV_\I)\|_\I$.}	
\end{figure}
$M=50$ and so $N=45$
Figure \ref{fig:3}(a) shows the influence of noise in the maximum norm.
Note that the solution $u^*\equiv 0$ corresponds to the finite-dimensional approximation $p^*\equiv0$,
which is a steady state of the ODE \eqref{eq:ODE}. The results show that the stochastic neighbourhood
defined by the maximum norm of the covariance matrix $\cV_\I$ widens as the bifurcation point $\mu=\mu_1^b$ 
is approached. This is due to critical slowing down, which is an effect well-studied in certain regimes near 
bifurcation points for SODEs \cite{KuehnCT2}. In particular, the increase in the (co-)variance can be used 
as an early-warning sign for the upcoming bifurcation point. For SPDEs, the scaling of the covariance 
operator $V_\I$ has been investigated recently using analytical techniques in one-dimensional normal 
form-type SPDEs for pattern formation near branch points \cite{GowdaKuehn}, where it is shown that we 
formally expect a scaling 
\be
\label{eq:earlywarning}
\langle V_\I e_k,e_j\rangle=\frac{\kappa}{\mu_1^b-\mu}=\cO\left(\frac{1}{\mu_1^b-\mu}\right),
\qquad \text{as $\mu\nearrow \mu_1^b$}
\ee
in certain modes $e_k$ and $e_j$, where $\kappa$ is a constant. Since we are taking the maximum norm here, 
the result \eqref{eq:earlywarning} does not apply directly but we still see an increase in $\|\cV_\I\|$ 
that shows a qualitatively similar in the inset in Figure \ref{fig:3}(a), where we plot \eqref{eq:earlywarning}
with $\kappa=0.001$. However, there are differences near the bifurcation point, which still have to be
explored in future work; see also \cite{GowdaKuehn} where this issue has been partially discussed. Of course,
very close to a bifurcation point, one has to take into account the full nonlinear dynamics of the problem
and the linearization approach via Ornstein-Uhlenbeck processes is no longer valid; we also leave this as
an extension for future work.

\begin{figure}[htbp]
\centering
\psfrag{V}{$u$}
\psfrag{x}{$x$}
\psfrag{y}{$y$}
		\includegraphics[width=0.75\textwidth]{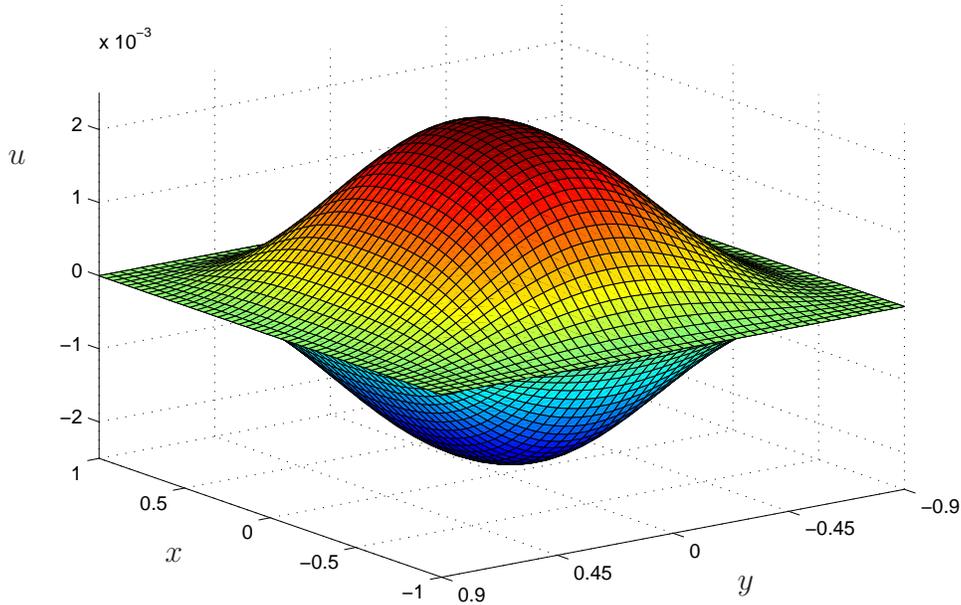}
\caption{\label{fig:5}Visualization of the stochastic neighbourhood calculated from
the covariance matrix for $\mu=1$ with $\varphi_k=k$ for $k\in\{1,2\}$. Each point in the 
grid corresponds to a point 
$p_j=u^{m,n}=u(x_m,x_n)$ and on the vertical axis labelled $u$, we show the stationary 
solution $u^*\equiv 0$ (a green plane) as well as the associated variances for each 
point $(\cV_\I)_{jj}$ (the surface with coloring green/yellow/red) and $-(\cV_\I)_{jj}$ 
(the surface with coloring green/cyan/blue). The mesh is overlayed on all three objects.
Note that it is difficult to visualize the full evolution of the associated variances 
when varying the parameter $\mu$ in a three-dimensional plot so this paper contains as
a supplementary material a movie, which displays this figure for $\mu\in[1,\mu_1^b)$.}	
\end{figure}

Figure \ref{fig:3}(b) shows a zoom near the covariance neighbourhood boundary curves defined by 
$\|\cV_\I\|_\I$ to illustrate that there is only a very small difference between the three noise 
truncation levels $K=2,4,8$. In fact, the differences are very close to the numerical precision used for the
linear solver $\texttt{tollin}=10^{-4}$. This result is expected since the largest stochastic
fluctuation, when measured in the $\|\cdot\|_\I$, is dominated by the mode with the largest eigenvalue.
The largest eigenvalue is included in all three cases $K=2,4,8$ as shown in Figure \ref{fig:3}(c).\medskip

Figure \ref{fig:4} compares the growth of the stochastic fluctuations in different norms. In this 
case, we focused only on the vector of variances 
\be
\text{diag}(V_\I):=\{(\cV_\I)_{jj}\}_{j=1}^J
\ee
Note carefully that the scale difference in Figure \ref{fig:4} is over three orders of magnitude, so 
one has to be very careful to define a threshold for covariance neighbourhoods only when it is clear, 
which norm is used to evaluate this threshold. The most important aspect of Figure \ref{fig:4} is that
in all the different norms, an early-warning sign, based upon the scaling law of the covariance neighbourhoods
before the bifurcation point, is clearly visible.\medskip

\begin{figure}[htbp]
\centering
\psfrag{T}{$T$}
\psfrag{mu}{$\mu$}
\psfrag{a}{(a)}
\psfrag{b}{(b)}
		\includegraphics[width=1\textwidth]{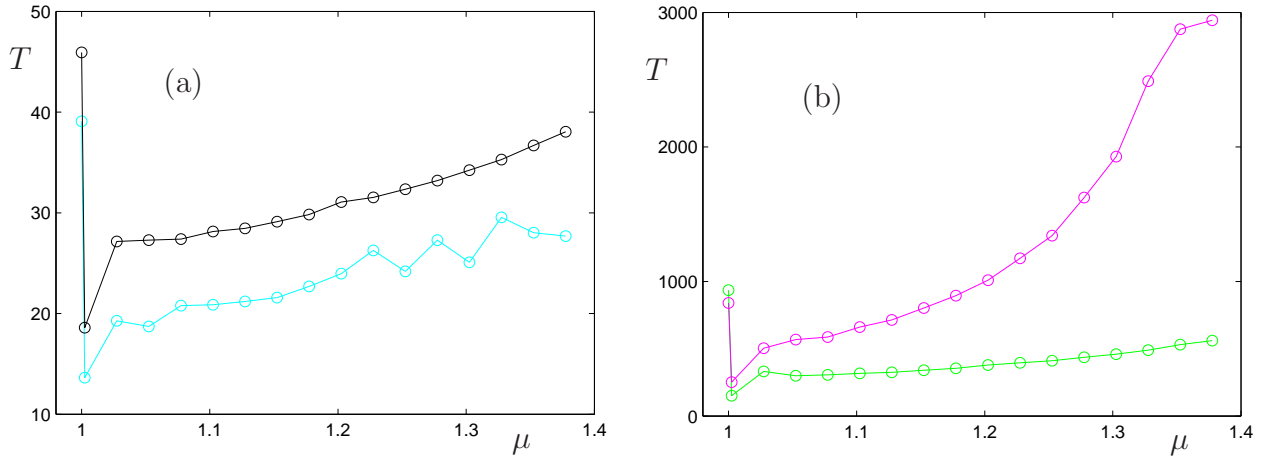}
\caption{\label{fig:9}Computation time $T$ for Lyapunov operator
computation for the lower homogeneous branch with $\mu\in[0,\mu_1^\txtb)$
for three different iterative linear solvers. 
Each circle represents a point on the continuation curve, where \eqref{eq:Lya1} has been 
solved using \texttt{bicgstab} (cyan), \texttt{gmres} (green/magenta) and \texttt{qmr} (black). 
For \texttt{gmres}, there is also the option to restart after a fixed number of inner iterations,
which has been set to 10 (magenta) and 0 (green) respectively. The linear interpolation between
the points is for visualization purposes. The parameter $\mu$ is plotted against computation time
$T$ (in seconds) in (a) and (b). Note that the \texttt{gmres} method in (b) has a lot longer 
computation time.}	
\end{figure}

Of course, one may also consider a visualization of the results in three dimensions. In Figure \ref{fig:5}
the stationary solution $u^*\equiv 0$ is shown together with a stochastic neighbourhood defined via the variances
$\text{diag}(\cV_\I)$. More precisely, each grid point $(x_m,y_n)\in \cD$ has an associated solution value
$u(x_m,y_n)=p_j$ and an associated variance $(\cV_\I)_{jj}$; of course, we plot also $-(\cV_\I)_{jj}$. One clearly
sees that the fluctuations have been computed for the Dirichlet Laplacian with zero boundary conditions
so that the fluctuations are largest at the center of the domain. Figure \ref{fig:5} shows the case when $\mu=1$;
in the supplementary material a movie is provided showing the growth of the stochastic neighbourhood in
three dimensions upon variation of $\mu\in[1,\mu_1^b)$; for a brief comparison to Monte-Carlo simulations of 
the full nonlinear system see Appendix \ref{ap:1}.\medskip

\begin{figure}[htbp]
\centering
\psfrag{T}{$T$}
\psfrag{iter}{$\cI$}
\psfrag{mu}{$\mu$}
\psfrag{u}{$\|\cdot\|$}
\psfrag{a}{(a)}
\psfrag{b}{(b)}
\psfrag{c}{(c)}
		\includegraphics[width=1\textwidth]{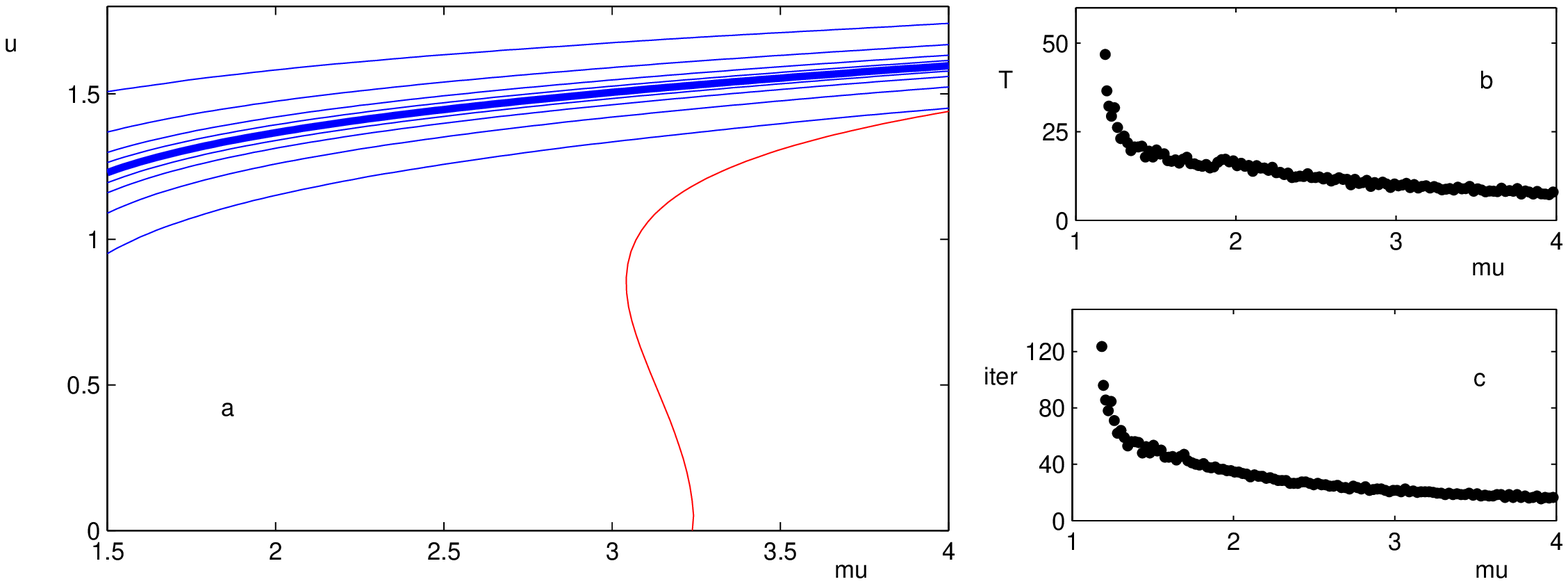}
\caption{\label{fig:6}Computation of the upper locally asymptotically stable branch for 
$\mu\in R_1\cup R_2\cup R_3\cup R_4$. (a) All norms in this figure are maximum norms $\|\cdot\|=\|\cdot\|_\I$. 
We do not show the initial part in region $R_1$ near the fold point in this part of the figure; see Figure 
\ref{fig:7}. The deterministic steady state solution $u^*$ on the branch $\Gamma_1$ (thick blue curve) 
from the deterministic continuation run in Figure \ref{fig:1} is shown together with the unstable branch 
$\Gamma_2$ (red). The branches $\Gamma_1$ and $\Gamma_2$ correspond to steady states $p^*$ of the discretized 
system \eqref{eq:ODE}. Furthermore, the norm of the covarince matrix $\cV_\I$ is plotted to indicate the size 
of the fluctuations as measured by the linearized system (thin blue curves). There are actually four different 
neighbourhoods shown for additive noise, eigenvalues \eqref{eq:noise_num} for 
$\tilde{\sigma}=200, 400, 800, 1600$. (b) Computation time $T$ along the entire upper branch. (c) Number of 
iterations $\cI$. In (b)-(c) each dot represents a point on the continuation curve,
where \eqref{eq:Lya1} has been solved using the \texttt{bicgstab}-method. Only the computation for 
$\tilde{\sigma}=200$ is shown, the other three cases yield very similar results.}	
\end{figure}

The performance of different standard iterative linear solvers for the large sparse linear 
system \eqref{eq:Lya_linsyst} has also been investigated. The lower homogeneous branch has 
again been used as a test case. We used \texttt{bicgstab}, \texttt{gmres} and \texttt{qmr} in their
MatLab implementation\footnote{Version: MatLab R2013a} with tolerance $10^{-4}$. As before, we use 
the eigenvalues \eqref{eq:noise_num} now with $\varphi_k=(k-1)\frac{4}{10}$ for $k\in \{1,2,\ldots,11=K\}$. 
Figure \ref{fig:9} shows a comparison in the computation time for each point on the bifurcation curve. 
The results show that there are certainly differences in the computation time for the different
methods and that the \texttt{bicgstab} method is the fastest in our context while 
\texttt{gmres} is the slowest for this problem. Given the different performances
observed here for various iterative linear system solvers for the Lyapunov equations we have to 
solve, we can conjecture that using a specialized algorithm for Lyapunov equations may generate 
substantial additional speed-up as indicated in (S1) above.

\subsection{The Upper Branch}
\label{ssec:upper_res}

In the next numerical continuation run, we focus on the upper locally asymptotically stable steady state
branch $u=u^*$ computed in Section \ref{sec:det_cont}; see Figure \ref{fig:1}. The algorithmic parameters
for the linear system solvers are kept as for the homogeneous branch computations from Figures 
\ref{fig:2}-\ref{fig:5} now with $\varphi_k=(k-1)\frac{8}{20}$ for $k\in \{1,\ldots,21=K\}$ 
and eigenvalues given by \eqref{eq:noise_num} for additive noise. Furthermore, we recall
that there is a locally asymptotically stable non-trivial solution branch $\Gamma_1$ in the parameter
regime $\mu\in R_1\cup R_2 \cup R_3\cup R_4$ as shown in Figure \ref{fig:1}. The branch $\Gamma_1$ changes 
stability at the fold point $\mu=\mu_1^f\approx 1.1798$.\medskip

\begin{figure}[htbp]
\centering
\psfrag{V}{\small $\|\cdot\|$}
\psfrag{mu}{$\mu$}
\psfrag{muc}{$\mu_1^f-\mu$}
\psfrag{u}{$\|\cdot\|$}
\psfrag{a}{(a)}
\psfrag{b}{(b)}
		\includegraphics[width=1\textwidth]{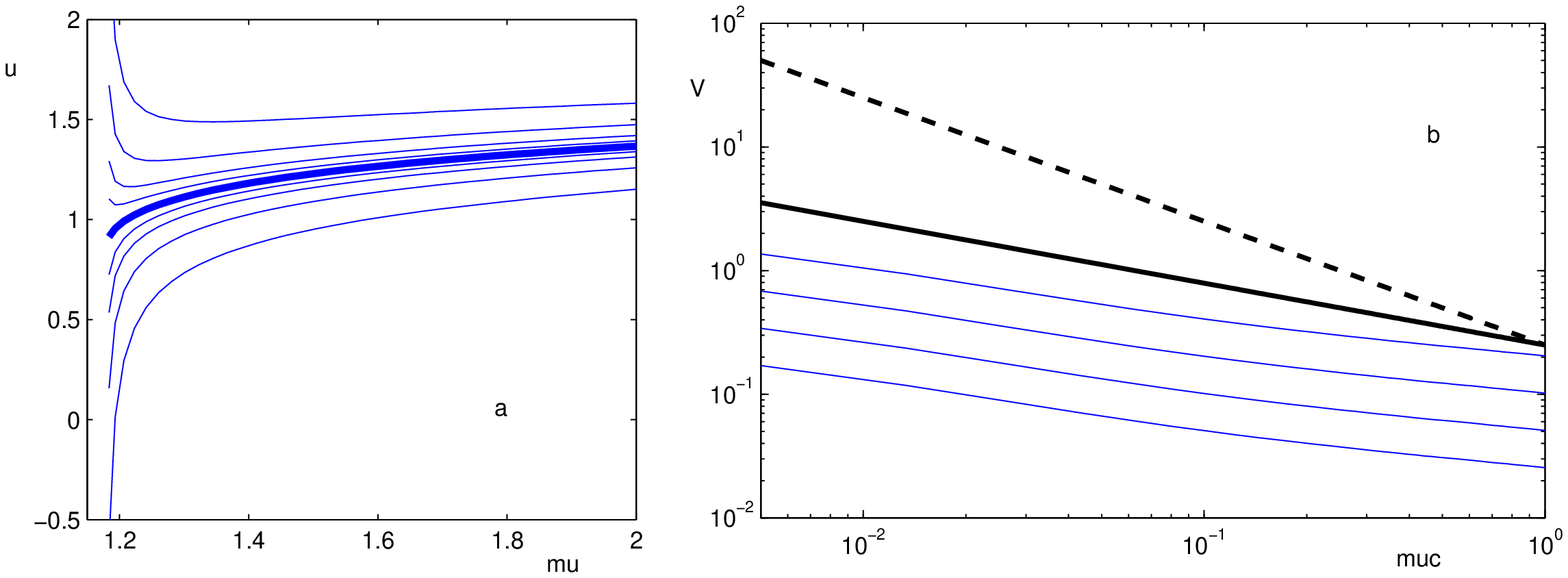}As before, we use 
the eigenvalues \eqref{eq:noise_num} now with $\varphi_k=(k-1)\frac{4}{10}$ for $k\in \{1,2,\ldots,11=K\}$. 
\caption{\label{fig:7}Computation of the upper locally asymptotically stable branch for 
$\mu\in R_1\cup R_2\cup R_3\cup R_4$, now with a focus on the region $R_1$ near the fold point 
at $\mu=\mu_1^f\approx 1.1798$. (a) All norms in this figure are maximum norms $\|\cdot\|=\|\cdot\|_\I$. The 
color coding of (a) is the same as in Figure \ref{fig:6}(a). There are four different 
neighbourhoods shown for additive noise, eigenvalues \eqref{eq:noise_num} for 
$\tilde{\sigma}=200, 400, 800, 1600$. (b) Computation of the scaling laws for the growth of the norm 
$\|\cV_\I\|_\I$ as the bifurcation point is approached. The four thin blue curves are just a 
double-logarithmic plot of the curves from (a) with the bifurcation point shifted to zero. The thick black 
dashed line has slope $-1$, while the thick black solid line has slope $-\frac12$; the latter scaling 
fits the observed scaling law very well.}	
\end{figure}

Figure \ref{fig:6}(a) shows the continuation run with the branch $\Gamma_1$ and a focus on 
$\mu\in R_2\cup R_3 \cup R_4$, {i.e.}, bounded away from the fold point. Four stochastic neighbourhoods 
computed for different noise strengths are shown. If we increase $\mu$, then the deterministic stability
increases. Since the additive noise is kept constant, this explains the shrinking of the four neighbourhoods.
Note carefully, that the neighbourhoods will interact with the unstable branch $\Gamma_2$.
The precise interaction threshold depends upon the nonlinearity but we can observe that there is a 
competition effect: the norm of the solution on the unstable branch increases while the size of the 
stochastic neighbourhoods decreases. In the current setup, we can observe that the noise-induced 
effect is going to dominate eventually for sufficiently large $\mu$ and we expect frequent noise-induced 
excursions on a sub-exponential time-scales across the unstable branch for individual sample paths. 

Figures \ref{fig:6}(b)-(c) show the computational performance along the branch $\Gamma_1$. As we start 
with the zero vector as an initial guess for the first value smallest value of $\mu$ near $\mu_1^f$, 
the iterative method is initially slow but again increases performance after a few points on the 
continuation branch as in Figure \ref{fig:2}. We also observe from the results that it seems to be 
computationally faster to compute the Lyapunov equation in a regime, where the leading eigenvalue of 
the linearized deterministic system is smaller, {i.e.}, in a more stable regime.\medskip

{\small \textbf{Remark:} In Figures \ref{fig:6}-\ref{fig:7} we have increased the noise level by several
orders of magnitude to visualize the neighbourhoods better in comparison to the order one scale
of the deterministic bifurcation branch. Note that the Lyapunov equation is invariant under a 
change of scale, i.e., a multiplicative factor $\sqrt{\zeta}$ for $\zeta>0$ in front of the noise 
term defined by $\cB$ enters as a scale factor of $1/\zeta$ for the covariance $\cV$ in; see the 
Lyapunov equation \eqref{eq:Lya_linsyst}. This avoids having to add several zoom-ins near the bifurcation
branch in Figures \ref{fig:6}-\ref{fig:7}.}\medskip 

Figure \ref{fig:7} focuses on the fold point. In Figure \ref{fig:7}(a), we consider $\mu\in R_1$
and illustrate the growth of the stochastic fluctuations near the bifurcation point. Again, we are 
interested in scaling laws. For example, the goal is to determine the exponent $\alpha>0$ in 
\be
\label{eq:earlywarning1}
\|\cV_\I\|_\I=\frac{\kappa}{(\mu-\mu_1^f)^\alpha}=\cO\left(\frac{1}{(\mu-\mu_1^f)^\alpha}\right),
\qquad \text{as $\mu\searrow \mu_1^f$.}
\ee
Figure \ref{fig:7}(b) shows that $\alpha=1$ does not yield a good fit near the fold in contrast to
the case near the bifurcation from the homogeneous branch at $\mu=\mu_1^b$. The fit for $\alpha=\frac12$
is a lot better and this corresponds precisely to the theory expected from SODEs. Indeed, near 
transcritical and pitchfork bifurcations one finds for SODEs $\alpha=1$, while for the fold one
obtains $\alpha=\frac12$ \cite{KuehnCT2}. Hence, the result in Figure \ref{fig:7}(b) is a natural
generalization from the SODEs results as expected from the recent SPDEs results on early-warning 
signs in \cite{GowdaKuehn}.\medskip

\begin{figure}[htbp]
\centering
\psfrag{mu}{$\mu$}
\psfrag{u}{$\|\cdot\|$}
\psfrag{a}{(a)}
\psfrag{b}{(b)}
\psfrag{c}{(c)}
		\includegraphics[width=1\textwidth]{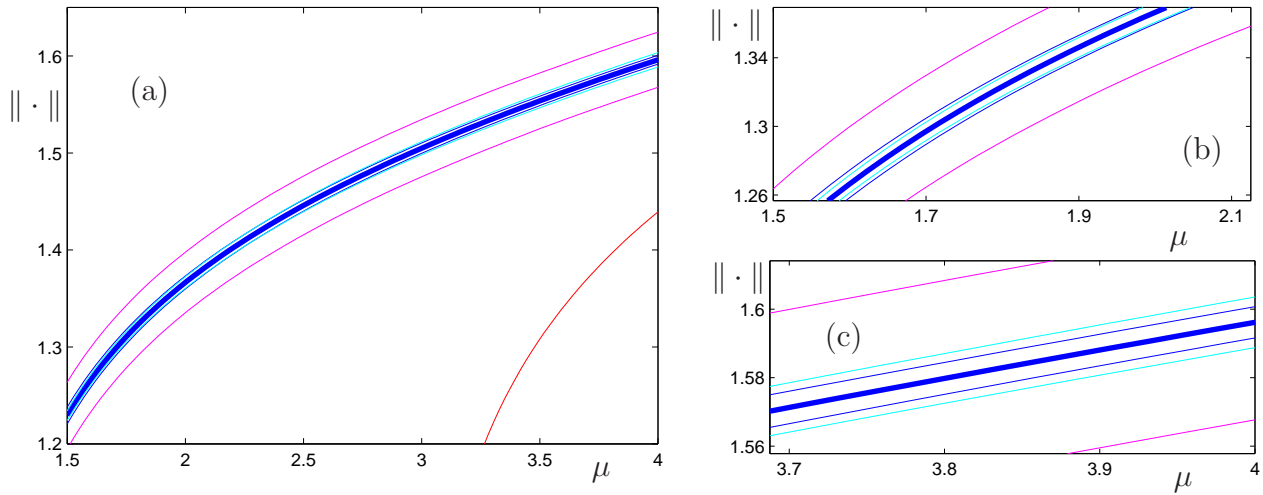}
\caption{\label{fig:8}Results for different types of scalar-valued noise terms for fluctuations around 
the upper locally asymptotically stable deterministic bifurcation branch with 
$\mu\in R_1\cup R_2\cup R_3\cup R_4$. All norms in this figure are maximum norms $\|\cdot\|=\|\cdot\|_\I$ 
(a) The deterministic homogeneous steady state solution $u^*$ on the branch $\Gamma_1$ (thick blue curve) 
from the deterministic continuation run in Figure \ref{fig:1} is shown together with the unstable 
branch $\Gamma_2$ (red). Three different noise terms are considered and the associated neighbourhoods
are plotted: additive noise $G(u)\equiv 1$ (thin blue curves), norm-dependent noise scaling 
\eqref{eq:noise3} (thin cyan curves) and norm-dependent noise amplification/reduction \eqref{eq:noise4} 
(thin magenta curves). The value of $\tilde{\sigma}=50$ has been fixed and the 
eigenvalues \eqref{eq:noise_num} with truncation $K=20$ have been used. (b)-(c) Zooms of different parts
of the curves in (a).}	
\end{figure}

In the next numerical continuation run we investigate different noises. The most natural first 
generalization from purely additive noise $G(u)\equiv 1$ is to consider scalar-valued noise terms 
$G(u)\in\R$, where $G$ does depend upon $u$, {i.e.}, relatively simple multiplicative noise terms. 
The choices we make here are purely for test purposes and many other situations may arise in
practical applications. We consider additive noise and the following two cases
\bea
G(u) & = & \frac12 \|u\|_\I^2, \label{eq:noise3}\\
G(u) & = & \|u\|_\I-u. \label{eq:noise4} 
\eea 
The term \eqref{eq:noise3} models a scaling effect of the noise depeding upon 
the current maximal value of the solution. The term \eqref{eq:noise4} is a shifted multiplicative 
noise term. The results 
in Figure \ref{fig:8} show that multiplicative noise can indeed
crucially change the relative behaviour in different regimes. For example, if $\mu\in R_1$ so that we 
are near the fold point, then the stochastic neighbourhood generated by the noise-amplification term
\eqref{eq:noise3} is smaller than the additive noise neighbourhood, while the situation is interchanged
for $\mu\in R_4$, {i.e.}, there is a change in the intermediate regimes when $\mu$ passes through 
$R_2\cup R_3$. Furthermore, the results show that it is a-priori not easy to see from the structure
of the noise term, which stochastic neighbourhood is larger when two different types of multiplicative
noise are compared. Hence, the method proposed here can automate this process, even as a fast post-processing
step for deterministic bifurcation diagrams.

\appendix

\section{Appendix: Approximation near Steady States}

In this entire work, we have always stressed that we work in the vicinity of deterministically
locally asymptotically stable steady states. The goal is to understand the influence of stochastic
forcing on the PDE in this regime before any large deviation effects have occurred, i.e., on
subexponential time scales. In this case, we have made the approximation that the full nonlinear 
SPDE can be approximated locally by the OU process \eqref{eq:SPDE_lin}. Although it is well-known
from previous numerical experiments from SODEs as well as from large deviation theory 
\cite{FreidlinWentzell}, and its versions and refinements \cite{BerglundGentz}, we provide for
completeness in this appendix also numerical tests for SPDEs and quote some relevant theoretical
results.

\subsection{Comparison to Monte-Carlo Simulation}
\label{ap:1}

In this section, we briefly compare the local dynamics obtained from the methods used in Section 
\ref{sec:contSPDE_Lya} with results from direct Monte-Carlo simulation. As a basic test case, 
we consider the same setup as in Section \ref{ssec:lower_res} focusing on the lower branch at $\mu=1$.
We use $\varphi_k=(k-1)\frac{4}{10}$ for $k\in \{1,2,\ldots,11=K\}$, $\tilde{\sigma}=5$, 
$(M,N)=(50,45)$, and the eigenfunctions \eqref{eq:efuncs}, as before. To
carry out the Monte-Carlo simulation, we integrate the full nonlinear system of SODEs derived 
in Section \ref{sec:fd} using a standard Euler-Maruyama method as described in \cite{Higham}.\medskip

\begin{figure}[htbp]
\centering
\psfrag{x}{$x$}
\psfrag{y}{$y$}
\psfrag{u}{$u$}
\psfrag{t1}{$t=0.25$}
\psfrag{t2}{$t=0.5$}
\psfrag{t3}{$t=0.75$}
\psfrag{t4}{$t=1$}
		\includegraphics[width=1\textwidth]{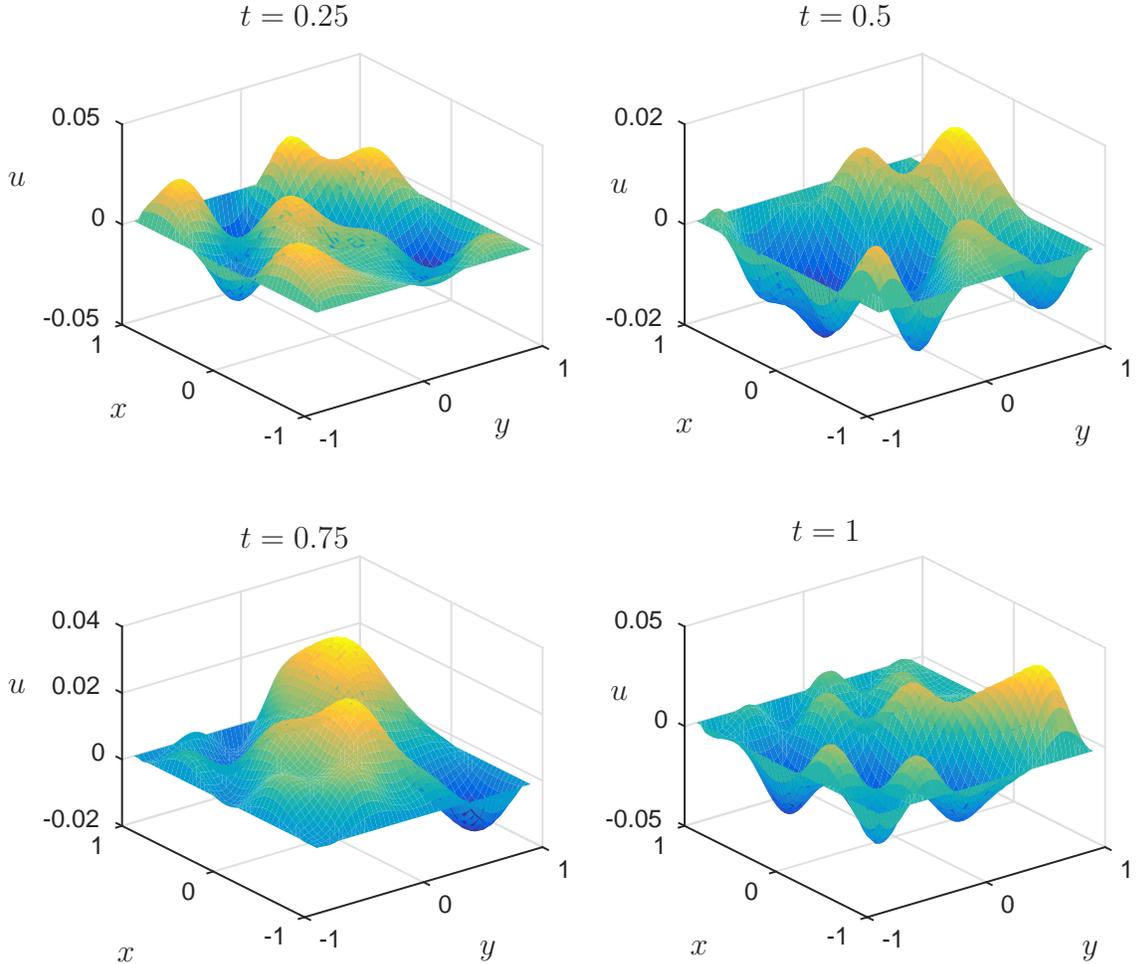}
\caption{\label{fig:10}Time snapshots of the solution obtained from Monte-Carlo simulation of the full nonlinear cqQC SPDE \eqref{eq:AC_SPDE} for parameters stated in the text in Section \ref{ap:1}.}	
\end{figure}

An output of a single sample path of the problem is shown in Figure \ref{fig:10} with snapshots at 
$t=0.25,0.5,0.75,1$. Due to the choice of the noise, the solution is relatively smooth and 
we also observe that there is no pattern-formation in the stable subthreshold regime visible 
for $\mu=1$. Basically, we only see small scale modulations due to the stochastic forcing.\medskip

\begin{figure}[htbp]
\centering
\psfrag{uemax}{$\cC_{\max}$}
\psfrag{t}{$t$}
\psfrag{u}{\scriptsize $\|\cdot\|_\I$}
\psfrag{uemin}{$\cC_{\min}$}
\psfrag{umin}{\bf \textcolor{green}{$\min$}$(u(x_m,y_n,t))$}
\psfrag{umax}{\bf \textcolor{red}{$\max$}$(u(x_m,y_n,t))$}
		\includegraphics[width=1\textwidth]{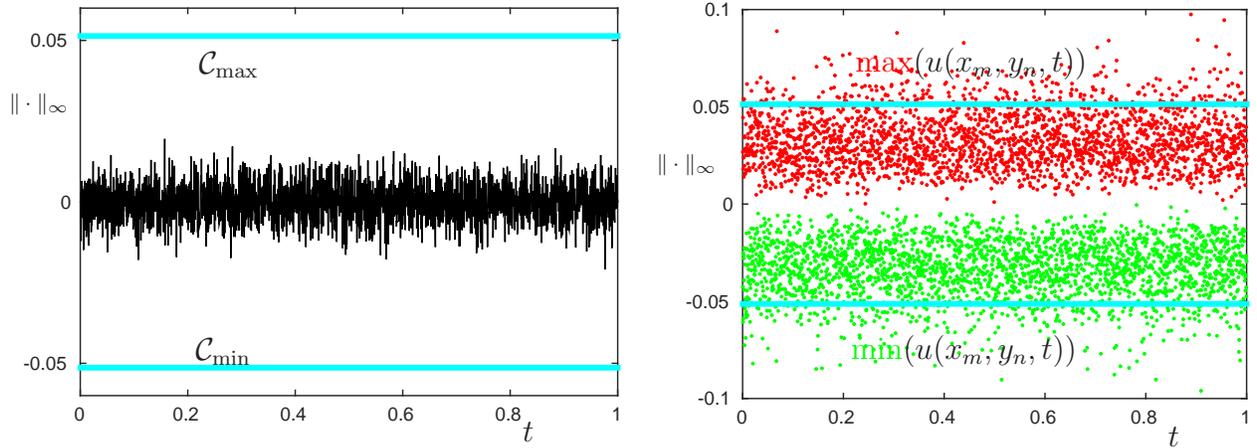}
\caption{\label{fig:11}Comparison between direct Monte-Carlo simulation and the 
neighbourhood $\cC(1)$ at $\mu=1$. The left figure shows the sample path $u(0,0,t)$ at the 
origin in the center of the domain $\cD=[-1,1]\times[-0.9,0.9]$. The cyan lines
are given by  $\cC_{\max}:=\max(\cV_{\I})$ and $\cC_{\min}:=\min(\cV_{\I})$ where the maximum
and minimum are calculated over all elements of the matrix $\cV_\I$ obtained by solving the 
Lyapunov equation. The right figure also shows the same cyan lines, now in comparison to 
maximum and minimum of the entire solution sample path over the domain as defined in \eqref{eq:maxmindom}, i.e.,
the red dots give the maxima at each time point and the green dots the minima.}
\end{figure}

To compare the results to the stochastic neighbourhoods defined by covariance operators, 
respectively the finite-dimensional approximation by matrices, consider again the 
definition \eqref{eq:Bh_big}. For $\mu=1$, we compute $\cV_\I$ as before solving
the Lyapunov equation \eqref{eq:Lya_linsyst} which yields $\cC(r)$. We fix $r=1$ as the 
radius and consider the $\|\cdot\|_\I$-norm. Figure \ref{fig:11} shows the results of this
comparison. We observe that the fluctuations of $u(x_m,y_n,t)$ at a generic single point in 
the domain are typically well-contained within $\cC(1)$ on a time scale of order $1$.
Furthermore, we have studied in Figure \ref{fig:11} the maximum and minimum of the 
Monte-Carlo simulation at the mesh points, i.e., 
\be
\label{eq:maxmindom}
\max_{m,n}u(x_m,y_n,t)\qquad \text{and} \qquad \min_{m,n}u(x_m,y_n,t). 
\ee
It is clearly visible that the maxima and minima are well-described by the boundaries 
computed from the neighbourhood $\cC(1)$. Indeed, on the time scale of order one only 
very few points $u(x_m,y_n,t)$ lie outside $\cC(1)$ while most points stay inside $\cC(1)$
for all times or just make very small excursions outside of $\cC(1)$. Only on time scales
substantially larger than order one, we expect that large deviations going far beyond $\cC(1)$
occur as discussed from a theoretical perspective in the next section.

\subsection{Remarks on Theoretical Approximation}
\label{ap:2}

From an abstract perspective, we expect the local approximations we made to work well on
certain time scales for small noise because of the fundamental results on large deviation 
theory \cite{FreidlinWentzell}. For example, consider again the SODE \eqref{eq:SODE} with additive noise and 
fixed parameter values given by 
\benn
\txtd p=\vartheta(p)~ \txtd t+\sigma~\txtd \beta,
\eenn 
where $p\in\R^J$, $\vartheta$ is a sufficiently smooth vector field, $\sigma\geq 0$ determines
the noise level and $\beta$ is a $J$-dimensional vector of iid Brownian motions. Let $\cK\subset \R^J$
be a bounded domain with smooth boundary and consider the case $p(0)\in \cK$. Define the first exit time by 
\benn
\tau^\sigma_{\cK} :=\inf\{t:p(t)\not\in \cK\}.
\eenn
On absolutely continuous paths $\gamma\in C([0,T],\R^J)$ define the action functional
\benn
S(\gamma)=\frac12 \int_0^T \|\dot{\gamma}(s)-\vartheta(\gamma(s))\|_2^2~\txtd s.
\eenn
Let $\cH=\cH(t,q)$ denote the space of all absolutely continuous paths $\gamma$ with $\gamma(0)=q$
and $\gamma(s)\not\in \cK$ for some $s\in[0,t]$. Then one has \cite[Ch.4,Thm.1.2]{FreidlinWentzell}
\be
\label{eq:ldFW}
\lim_{\sigma\ra 0} \sigma^2\ln \P(\tau^\sigma_{\cK}\leq t| p(0)=q)=-\min_{\gamma\in\cH}S(\gamma).
\ee  
The main message from \eqref{eq:ldFW} is that the expected escape time from a suitable neighbourhood 
$\cK$ of a deterministically locally asymptotically stable steady state scales,  exponentially 
as $\cO(\exp(-K/\sigma^2))$ in the small noise limit for some constant $K>0$. Hence, it will usually take very long until 
a path $p$ leaves the ellipsoidal neighbourhoods we used, e.g., just pick $\cK=\cC(r)$ for a suitable 
radius $r>0$. For long periods we expect the linear approximation via an OU process to hold near a 
deterministically locally asymptotically stable steady state. This idea has been used to develop
refinements of large deviation ideas. In particular, there are known
results that have the structural form \cite[Thm.5.1.6]{BerglundGentz}
\be
\label{eq:BG}
\P(\tau_{\cC(r)}\leq t|p(0)\in\cC(r))\leq K_1\exp(-K_2 r^2/\sigma^2)
\ee   
where $K_1=K_1(r,\sigma,J,t)$, $K_2=K_2(r,\sigma)$ are functions that depend upon the given problem, i.e.,
on noise level $0<\sigma\ll1$, the scaling of the neighbourhood $\cC(r)$, and the dimension $J$. The prefactor
$K_1$ obviously does grow as $t$ increases but for $K_{1,2}$ explicit estimates are available. Basically, a 
result of the form \eqref{eq:BG} allows us to interpret the neighbourhoods $\cC(r)$ as confidence 
neighbourhoods for the full nonlinear problem up to a certain large time. Now the remaining question is
whether the ideas available for SODEs can also be transferred to the setting of SPDEs. However, as long as
we use a numerical scheme to convert an SPDE into a large set of $J$ SODEs first, it is natural to expect that
the additional error in estimates of the form \eqref{eq:ldFW}-\eqref{eq:BG} is a numerical approximation
error that can be controlled as $J\ra \I$; although there are no general schemes available to make this 
idea precise, recent initial findings are encouraging \cite{BloemkerJentzen,BerglundGentz10}. Therefore, it
seems reasonable to propose a practical numerical scheme that aims to calculate $\cC(r)$ for wide classes 
of SPDEs as we have done in this work.

\bibliographystyle{plain}
\bibliography{../my_refs}

\end{document}